\documentclass[conference]{IEEEtran}
\IEEEoverridecommandlockouts
\usepackage{cite}
\usepackage{amsmath,amssymb,amsfonts}
\usepackage{algorithmic}
\usepackage{graphicx}
\usepackage{wrapfig}
\usepackage{textcomp}
\usepackage{xcolor}

\begin{document}

\newtheorem{theorem}{Theorem}[section]
\newtheorem{corollary}{Corollary}
\newtheorem{lemma}{Lemma}
\newtheorem{definition}{Definition}[section]
\newtheorem{observation}{Observation}
\newtheorem{notation}{Notation}
\newtheorem{example}{Example}[section]	
\newtheorem{remark}{Remark}

\newcommand{\cross}[1]{%
\mbox{\vbox{\kern 1pt\hbox{\vbox{\hrule
\kern 2pt\hbox{\ensuremath{\vphantom{b}#1}\kern 2pt}}\vrule\kern 1pt}}}\,}

\title{The BF Calculus and the Square Root of Negation\\
\thanks{Kauffman's work was supported by the Laboratory of Topology and Dynamics, Novosibirsk State University (contract no. 14.Y26.31.0025 with the Ministry of Education and Science of the Russian Federation).}
}

\author{\IEEEauthorblockN{Louis H. Kauffman}
\IEEEauthorblockA{\textit{Department of Mathematics, Statistics and Computer Science} \\
\textit{University of Illinois at Chicago}\\
Chicago, IL USA \\
\textit{Department of Mechanics and Mathematics} \\
\textit{Novosibirsk State University, Novosibirsk, Russia}\\
kauffman@uic.edu
\and
\IEEEauthorblockN{Arthur M. Collings}
\IEEEauthorblockA{\textit{Independent Researcher} \\P.O. Box 114\\
Red Hook, NY USA \\
otter@mac.com}
}

%\and
%\and
%\IEEEauthorblockN{4\textsuperscript{rd} Given Name Surname}
%\and
%\IEEEauthorblockN{5\textsuperscript{rd} Given Name Surname}
%\and
%\IEEEauthorblockN{6\textsuperscript{rd} Given Name Surname}
}

\maketitle

\begin{abstract}
The concept of imaginary logical values was introduced by Spencer-Brown in Laws of Form, in analogy to the square root of -1 in the complex numbers. In this paper, we develop a new approach to representing imaginary values. The resulting system, which we call BF, is a four-valued generalization of Laws of Form. Imaginary values in BF act as cyclic four-valued operators. The central characteristic of BF is its capacity to portray imaginary values as both values and as operators. We show that the BF algebra   is an stronger, axiomatically complete extension to Laws of Form capable of representing other four-valued systems, including the Kauffman/Varela Waveform Algebra and Belnap\textquotesingle s Four-Valued Bilattice. We conclude by showing a representation of imaginary values based on the Artin braid group.  

\end{abstract}


\section{Introduction}\label{sec1}

We reexamine the concept of \emph{imaginary logical values}, first proposed by George Spencer-Brown \cite{bf1} and later re-conceptualized by Kauffman and Varela  \cite{bf3, bf4, bf5}. Imaginary logical values are analogous to $i = \sqrt{-1}$ in ordinary algebra. We introduce a new approach to imaginary Boolean values based on the concept of the \emph{square root of negation} as introduced by Kauffman in 1989 \cite{bf5, bf6} and revived more recently as a logical calculus that we call BF \cite{bf7, bf8, bf81}. 

The reader may recall that the square root of minus one, i, can be represented as acting on ordered pairs of real numbers by the formula $i(a,b) = (-b,a)$. This suggests a corresponding ``logical operator" of the form $\cross{(a,b)}_i = (\sim b,a)$ where $a$ and $b$ are elements of a Boolean algebra. This is one of the key ideas behind our construction of BF in this paper, which we stress at this point to indicate the close analogy of our work with the extension of the real numbers to the complex numbers.

The first section of the paper begins with a brief introduction to \emph{Laws of Form}, to the idea of distinction, and to the notation used by Spencer-Brown. We have emphasized the idea that values in Laws of Form can be indicated both as names and as operations. In Section III, we introduce several known four-valued calculi, which the reader may regard as precursors to BF. Section IV introduces the BF Calculus, stressing the following important points: 

\begin{itemize}
    \item Every consequence in Spencer-Brown\textquotesingle s Primary Algebra can be represented and formally demonstrated  in BF.
    \item Consequences exist and can be demonstrated in BF that do no exist in the Primary Algebra.
    \item All values (real and imaginary) can be viewed as being both a value and an operation. 
    \item Imaginary values act as cyclic operators, in exact analogy to a clock marked at the four quarter hours. 
\end{itemize}

Section IV concludes with a proof that the BF algebra is axiomaically complete based on an extended set of the Bricken Axioms. Finally, in Section V, we show that each of the four four valued calculi can be represented in BF, and that the BF Algebra is equivalent to a Four Valued Billatice.

\section{Laws of Form}

\subsection{Distinction}

\emph{Laws of Form} by Spencer-Brown \cite{bf1}, and the Primary Algebra (PA) it describes, is based on the idea of \emph{distinction}, represented by the dividing of a space into two regions, one  \emph{marked}, the second \emph{unmarked}. In \emph{Laws of Form}, the mark \cross{\quad} indicates the \emph{marked state}, and the empty value \quad \quad indicates the \emph{unmarked state}. The step of representing a value by an empty space, by the lack of a sign, is motivated by a key idea: doing so permits the  mark \cross{\quad} to act both as the name of a value and as an operation.

\begin{figure}[htbp]
\centerline{ \includegraphics[width=1.5in]{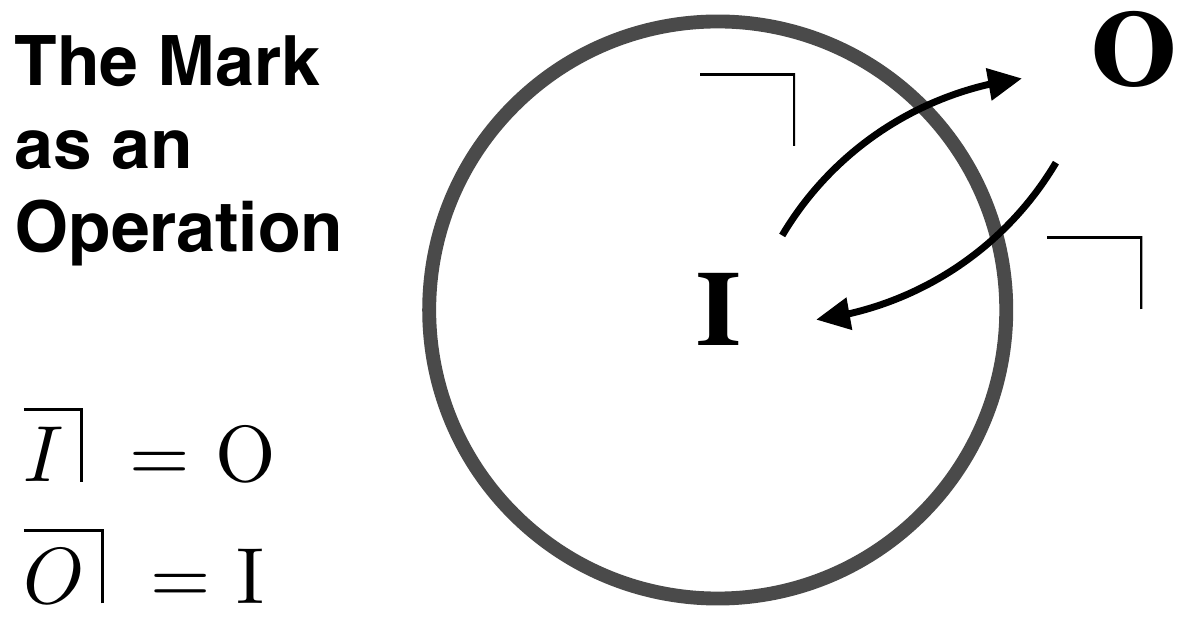}}
\caption{Representing a Distinction between Inside (I) and Outside (O)}
\label{fig1}
\end{figure}

Consider Figure \ref{fig1}, in which we have drawn a closed circle,  creating a distinction between inside, $I$, and outside $O$. We regard the mark \cross{\quad} as an operator that takes $I$ to $O$ and $O$ to $I$. Then we observe the following:
\begin{equation}\label{eq1}
\cross{I} = O, \quad \quad \cross{O} = I,
\end{equation}
\begin{equation}\label{eq2}
 \cross{\cross{I}} = \cross{O} = I, \quad \quad \cross{\cross{O}} = \cross{I} = O,
\end{equation}

\noindent 
 so for any state $X$ we have $\cross{\cross{X}} = X$.\\[-9pt]

The conceptual shift is to designate the inside to be unmarked (literally to have no symbol), so that $I =   \quad \quad  $. Then from (\ref{eq1}) we obtain
\begin{equation}\label{eq3}
 \cross{\quad} = O, \quad \quad \cross{O} = \quad \quad ,
\end{equation}
which means we have equated the mark \cross{\quad} with the outside $O$. From (\ref{eq2}), we obtain 

\begin{equation}\label{eq4}
 \cross{\cross{\quad}} =  \quad \quad. 
\end{equation}
By identifying the value of the outside with the result of crossing from the unmarked inside, Spencer-Brown has introduced a multiplicity meanings to the mark. The statement $\cross{\quad} = \cross{\quad}$ can be interpreted on the left side to mean ``cross from the inside" and on the right  as ``the name of the outside". 

The mark itself can be seen to divide its surrounding space into an inside and an outside. When we write $\cross{\quad} \cross{\quad}$, the two marks are  positioned mutually outside each other, and we can choose to interpret either mark as a name that refers to the outside of the other. We may also  interpret two such juxtaposed marks to indicate successive naming of the state indicated by the mark. In either case we can take \cross{\quad}\ \cross{\quad}  as an instance of the principle that to repeat a name can be identified with a single calling of the name: 
\begin{equation}\label{eq5}
\cross{\quad}\ \cross{\quad} = \cross{\quad}.
\end{equation}

At this point we have a single sign \cross{\quad} representing both the operation of crossing the boundary of a distinction and representing the name of the outside
of that distinction. Furthermore, since the mark itself can be seen to make a distinction in its own space, the mark can be regarded as referent to itself and to the (outer side) of 
the distinction that it makes. The two equations  (\ref{eq4}) and (\ref{eq5}) represent these aspects of understanding a distinction and the signs that can represent this distinction. We will now see that the two equations and a natural formalism for expressions in the mark become a formal system that can be seen as an `arithmetic' for Boolean algebra.

\subsection{The Primary Arithmetic}

On the basis of these considerations, Spencer-Brown defines a very simple calculus, which he calls the \emph{Primary Arithmetic}. \\[-6pt]

\begin{definition}\label{def100}
An expression in the Primary Arithmetic is a finite pattern containing no variables satisfying

\begin{enumerate}
\item The empty pattern \quad \quad is an expression. 
\item  If $X$ is an expression, then $\cross{X}$ is an expression.
\item If $X$ and $Y$ are expressions, then $X\ Y$ is an expression. 
 \\[-6pt]
\end{enumerate}
\end{definition}

\noindent
Thus, an expression is a pattern of marks drawn in the plane such that of any two marks in the pattern one can say that neither is inside the other, or one of the two is inside the other one.\\[-6pt]

\begin{example}For example, the following are all expressions: \\[-6pt]
\begin{description}
\item  $\cross{\cross{\quad}}$ and $\cross{\cross{\cross{\quad}}}$\\[-6pt]
\item $\cross{\quad}$, \ and  $\cross{\quad}\ \cross{\quad}\ \cross{\quad}$\\[6pt]
$\cross{\cross{\cross{\quad}\ \cross{\quad}}\ \cross{\quad}\ \cross{\quad}}$

\end{description}
\end{example}

Calculating in the Primary Arithmetic is based on two rules, derived from equations (4) and (5) in the prior section, which may be used to build or simplify expressions.  \\[-6pt]

A1. $\cross{\quad}\ \cross{\quad} = \cross{\quad}$.\\[-6pt]

A2. $\cross{\cross{\quad}} =  \quad \quad$.  \\[-6pt]

Note that by making A1 and A2 into possible steps that transform one expression into another, we have begun the possibility of calculation. In mathematics we take for granted the idea of calculation, but it is exactly in situations where there are possible series of transforming steps that calculation arises.\\[-6pt]

\begin{definition}\label{def210} 
Two expressions in the Arithmetic are \emph{equal} if one can be obtained from the other by applying A1 and A2 in a finite series of steps.\\[-6pt]
\end{definition}

\begin{definition}\label{def200} 
The two \emph{simple expressions} are     the empty expression \quad \quad and the mark \cross{\quad}.\\[-9pt]
\end{definition}

The primary arithmetic proves extremely useful, and provides a playground, a sandbox for exploring the patterns and forms that arise from juxtaposing and nesting expressions based on these definitions. Consider the expressions\\[-6pt]

 \cross{\cross{\quad }\ \cross{\cross{\cross{\quad}}}\   },\quad \cross{\cross{\cross{\quad}}},\quad and \cross{\quad 
     	\cross{\cross{\quad} \cross{\quad}\ \cross{\quad}\ }\quad \cross{\cross{\cross{\quad}}}}\quad.\\[-6pt]

\noindent 
Applying rules A1 and A2 we can that see the first expression simplifies to the unmarked state:\\[-6pt]

\noindent
\cross{\cross{\quad }\ \cross{\cross{\cross{\quad}}}\   } \\[5pt]
\begin{tabular}{@{\hspace{5ex}}l l }
= \cross{\cross{\quad }\ \cross{\quad}}  &[A2]\\[3pt]
= $\cross{\cross{\quad}}$ &[A1]\\[3pt]
= $\quad \quad \quad$. &[A2]\\[10pt]
\end{tabular}\\[-6pt]

\noindent 
We encourage the reader to explore these forms, to invent new ones, and  to try their hand playfully at applying the rules. \\[-6pt]
	
\begin{theorem}\label{th100}(Spencer-Brown). Every arithmetic expression simplifies to one of the two simple expressions by application of rules $A1$ and $A2$. This simplification is unique.

\emph{Proof:} See  \cite{bf1}. \\[-6pt]
\end{theorem}

In fact, one can show that if two arithmetic expressions $X$ and $Y$ are obtained, one from another, by any sequence of steps $A1$ and $A2$ (simplifying or complexifying), then they both simplify to the same unique value. Thus the equivalence relation generated by $A1$ and $A2$ has two equivalence classes corresponding to the two simple expressions
 \quad \quad and \cross{\quad}.  Every arithmetic expression must belong to one or the other class, which we call {\it unmarked} and {\it marked}.

\subsection{Patterns in the Arithmetic}

\begin{theorem}\label{th-arithmetic}
The following equalities are always true for arbitrary arithmetic expressions $A$ and $B$:
\begin{enumerate}
    \item $A\ \cross{\quad} = \cross{\quad}$
    \item $\cross{\cross{A}} = A$
    \item $\cross{A}\ B = \cross{AB}\ B$
\end{enumerate}

\emph{Proof:} By Theorem \ref{th100} every arithmetic expression $A$ or $B$ simplifies to $\cross{\quad}$ or $\quad \quad $, so it is sufficient to consider cases where $A$ and $B$ are limited to these two values. Then in expression 1), if $A$ = \quad \quad, then $A\  \cross{\quad} = \cross{\quad}$. And if $A$ = \cross{\quad}, then $A\ \cross{\quad} = \cross{\quad}\ \cross{\quad} = \cross{\quad}$ by rule A1, which completes the proof for 1). Proofs for 2) and 3) are similar and are left for the reader.
\end{theorem}

\subsection{The Primary Algebra} 

Observing patterns in the primary arithmetic provides the basis to establish a new calculus that incorporates variables.  \\[-6pt]

\begin{definition}\label{exp_d1}
Algebraic expressions in the \emph{Primary Algebra (PA)} satisfy:  
\begin{enumerate}
    \item the empty value $\quad \quad$ is an expression, 
    \item Variables, $A, B, C, \dots   $ and $A_1, A_2, A_3 $ are expressions, 
    \item if $X$ and $Y$ are expressions then so is $XY$,
    \item if $X$ is an expression, the so is $\cross{X}$. \\[-6pt]
\end{enumerate}
\end{definition}

\begin{example}The following are each algebraic expressions: \\[-6pt]
\begin{description}
\item  $\cross{\cross{A}}$ and $\cross{\cross{\cross{A}B}C}$\\[-6pt]
\item $A\ \cross{\quad}$, \ and  $\cross{A}\ B $\\[6pt]
$\cross{\cross{\cross{A}\ \cross{BC}}\ \cross{D}\  \cross{EF}}$\\[-6pt]
\end{description}
\end{example}

\begin{definition}\label{def500} \emph{Principle of Substitution}: Let $E_{1}$ = $E_{2}$ be equivalent algebraic expressions with a variable $X$, and let $S$ be any expression. If we replace all occurrences of $X$ in $E_1$ and $E_2$ with $S$ to obtain $F_1$ and $F_2$, then $F_1 = F_{2}$ by substitution.\\[-6pt]
\end{definition}

\begin{remark}
A variable is a token or sign that stands for the presence or absence of an arithmetic expression. If we substitute any arithmetic expression for a variable $X$ in an algebraic expression, the result is another arithmetic or algebraic expression. In addition, once we have this idea, a variable $X$ in an algebraic expression can stand for another algebraic expression!  Thus we can regard a variable in an algebraic expression as standing for the presence or absence of an algebraic expression and allow the operation of substitution as described above.\\[-6pt]
\end{remark}

\begin{remark}
Each axiom is a given algebraic identity and hence is subject to the rule of substitution. This means that each axiom stands for a an infinite list of specific algebraic or arithmetic statements. For example $\cross{\quad}\ A = \cross{\quad}$ implies that $\cross{\quad}\ \cross{\quad} = \cross{\quad}$ and indeed that a mark placed next to any expression whatever can be replaced by a single mark. Note that, by substitution in B1 and B2, the Bricken axioms imply the  arithmetic rules A1 and A2 that we discussed above. For B1 replace the variable A by \cross{\quad} and for B2 replace the variable A by the empty expression.\\
\end{remark}

We will adopt the Bricken Axioms \cite{bf9} B1, B2, and B3 as an axiomatically complete basis for the Primary Algebra. \\

\noindent 
$B1. \quad  \cross{\quad} \ A = \cross{\quad}$\quad \quad \quad \ (Integration) \\[3pt]
$B2. \quad \cross{\cross{A}} = A$\quad \quad \quad \quad \quad (Reflexion)\\[3pt]
$B3.  \quad \cross{A}\ B =  \cross{A\ B}\ B$\quad \   (Generation)\\[-6pt]

\noindent 
(Following Spencer-Brown, we also assume $AB = BA$ implicitly rather than as an axiom.)\\

\begin{definition}\label{def400} \emph{Demonstration}. Let $E_{1}$ be an expression. If $E_{1}$ can be transformed into a second expression $E_{2}$ by applying B1, B2, or B3 in a finite sequence, then  $E_{1}$ = $E_{2}$. Such a set of steps is called a \emph{Demonstration} and the resulting equivalence is called a \emph{Consequence}. \\[-6pt]
\end{definition}

Consequences B4-B10 can be demonstrated using axioms B1-B3. We give the first and leave the rest for the reader. \\

\noindent
B4. $\cross{\cross{ \ A \ } \ A}$ =  \quad \quad. \quad \quad \quad \quad \quad  \quad \quad (Position)\\[5pt]
\begin{tabular}{@{\hspace{5ex}}l l }
= $\cross{\cross{\quad} \ A}$ &[B3, Generation]\\[3pt]
= $\cross{\cross{\quad}}$ &[B1, Integration]\\[3pt]
= $\quad \quad \quad$. &[B2, Reflexion]\\[10pt]
\end{tabular}\\

\noindent
B5.  $\cross{ \cross{\ A} \cross{\ B}  } \ C$ =  $\cross{ \cross{\ A \ C} \cross{\ B \ C}  }$\quad \quad (Transposition)\\[5pt]
B6. $\cross{\cross{\ A} B}  \ A  = A$.  \quad \quad \quad \quad \quad  \quad \quad \ (Occultation)\\[5pt]
B7. $A \ A = A$    \quad \quad \quad \quad \quad  \quad \quad \quad \quad \ \ (Iteration) \\[5pt]
B8. $\cross{ \cross{\ A} \cross{\ B}   }  \cross{\cross{\ A}  B } = A$.    \quad \quad   \quad \quad \ \  (Extension)\\[10 pt]	
B9. $\cross{\cross{ \cross{\ A} \ B} \ C} $ = $ \cross{\ A \ C} \cross{\cross{\ B} \ C}$.   \quad \ (Echelon)\\[5pt]	
B10. $\cross{\cross{A \  C} \ \cross{\ B \  \cross{\ C}} \ }$ = $ \cross{\cross{\ A} C} \cross{\cross{\ B} \cross{\ C} }$ (Cross Trans). \\[-6pt]

\begin{definition}\label{def600}
A \emph{proof} that two expressions $F$ and $G$ are equivalent consists of a valid logical argument that $F=G$. We say that expressions $F$ and $G$ are {\it arithmetically equivalent} if they take the same values (via arithmetic simplification) for all common choices of marked or unmarked states for their variables. Arithmetical equivalence is the same as equivalence via truth tables where the table of values is taken to be the marked and unmarked states rather than $true$ and $false$.\\[-6pt]
\end{definition}

\begin{theorem}\label{th200} (Spencer-Brown) The Primary Algebra (PA) is axiomatically complete: if $F$ and $G$ are arithmetically  equivalent algebraic expressions, then $F = G$ can be demonstrated using the axioms B1, B2, B3. (Spencer-Brown used the set B4, B5.) 

\emph{Proof:} See  \cite{bf1}.  
\end{theorem}

\subsection{Interpreting Logical Operations in the Primary Algebra}
The PA is algebraic, with derivations that are expressed as equivalence rather than by implication and inference by the rule of modus ponens. However, it is entirely possible to interpret the PA as a Propositional Logic, as follows: \\[-6pt]

\begin{definition}\label{def650} Interpretation for Propositional Logic\\[-6pt]

\noindent
Logic $\quad \quad \quad \; \longrightarrow \quad \quad  $   Primary Algebra \\[6pt]
  $A$ \emph{Or} $B \quad \  \quad \longrightarrow  \quad \quad   A\ B  \quad \quad$\\[3pt]
 $A$ \emph{And} $B  \quad \ \ \longrightarrow  \quad \quad   \cross{\cross{A}\ \cross{B}} $\\[3pt]
 \emph{Not} $A    \quad \quad \quad   \longrightarrow  \quad \quad   \cross{A} $\\[3pt]
 $A$ \emph{Implies} $B \ \longrightarrow  \quad \quad   \cross{A}\ B $\\[-6pt]
 \end{definition}

In this interpretation, the mark \cross{\quad} represents $true$ and the absence of a mark \quad \quad represents $false$. An algebraic expression that is arithmetically equivalent to
 $\cross{\quad}$ is said to be a {\it tautology}. This is exactly the same as the usual notion that a tautology has a truth table of only {\it true} values. The Completeness Theorem~\ref{th200} assures us that the consequences of the axioms B1, B2, B3 include all tautologies.  
 
\section{Four Valued Calculi }\label{sec3}

The concept of generating four-valued logics using pair-based notation is well known. We discuss several ways of treating such pairs before introducing the BF calculus. In each, values take the form $V = (v_1, v_2)$, where $v_1$ and $v_2$ may take on the values of the Primary Arithmetic, \quad \quad and \cross{\quad}. The four constant values in each are (\quad, \quad), (\cross{\quad}, \cross{\quad}), (\cross{\quad}, \quad), and (\quad, \cross{\quad}). Each four-valued system shares the corresponding juxtaposition operation in the Primary Algebra. \\[-6pt]

\begin{definition}\label{7300}The \emph{juxtaposition} $A B$ of two expressions $A$ and $B$ is defined $A B = (a_1,a_2)\ (b_1,b_2) = (a_1 b_1, a_2 b_2)$, where $a_1b_1$ and $a_2b_2$ are the corresponding juxtaposition operations as defined in the 2-valued Primary Algebra.\\[-6pt]  
\end{definition}

What distinguishes the systems are variations in the \emph{enclosure} or \emph{nesting} operation \cross{X}. 
\subsection{The Direct Product: System $PA \times PA$\\[-6pt]}
%\noindent
Enclosure is defined: $\cross{A} = \cross{(a_1, a_2)} = (\cross{a_1},\ \cross{a_2})$. \\[-6pt]

Each axiom and consequence B1, B2, and B3 is valid in $PA \times PA$. Every pair of expressions that are equivalent in the Primary Algebra are also equivalent in $PA \times PA$, in particular B3 and B4 above. This system can be interpreted as a 4-valued Boolean algebra.  

\subsection{The Kauffman/Varela DeMorgan Algebra}

Using the definition $\cross{A}_w=\cross{(a_1, a_2)}_w = (\cross{a_2}, \cross{a_1})$, Kauffman and Varela \cite{bf3} define $WF$, a DeMorgan Algebra in which neither the Law of the Excluded Middle (B4) nor the axiom B3 are valid. Based on earlier work presented by Kauffman at ISMVL in 1978 \cite{bf4}, the authors establish $WF$ is axiomatically complete using  alternative axioms:  \\[-6pt]

\noindent
W1. $\cross{ \cross{\ A}_w \cross{\ B}_w  }_w \ C$ =  $\cross{ \cross{\ A \ C}_w \cross{\ B \ C}_w  }_w$\quad  (Transposition)\\[6pt]
W2. $\cross{\cross{\ A}_w B}_w  \ A  = A$.  \quad \quad \quad \quad \quad  \quad \quad  \quad \quad \ (Occultation)\\[-6pt] 

Kauffman and Varela interpret the value $M = (\cross{\quad}, \cross{\quad}) = \cross{\quad}$ to be \emph{marked}, and the value $U = (\quad, \quad)$ to be \emph{unmarked}.  They interpret the values $I = (\cross{\quad}, \quad)$ and $J = (\quad, \cross{\quad})$ to represent opposite phases of the oscillations $M, U, M, U  \dots $ and $U, M, U, M \dots $, which they use to describe fixed points and waveform patterns resulting from \emph{reentering} and \emph{self-referential} equations. \\[-6pt]

Note that in the actual Kauffman/Varela algebra they use the notation $\cross{A}_w = \cross{A}$ with no extra designation. In fact, if one identifies the unmarked state with its 
empty parenthetical representative, $\quad \quad = (\quad, \quad),$ then we would have 
$$\cross{\quad \quad}_w = \cross{(\quad, \quad)}_w=(\cross{\quad}, \cross{\quad}) = \cross{\quad}$$
so that the new operator is identified notationally as an extension of the original operator. Since the new operator is in fact such an extension, this notation for it works well.

\subsection{Belnap\textquotesingle s DeMorgan Algebra}
 
Alternatively, the definition $\cross{A}_b = \cross{(a_1, a_2)}_b = (a_2, a_1)$ results in a second DeMorgan Algebra, Belnap\textquotesingle s $FOUR$, when combined with an alternative form of disjunction $A \vee B$. (See Section \ref{sec_v}.B.) Belnap interprets the value $(\quad, \cross{\quad}) = T$ as {\it true} (not {\it false}), the value $(\cross{\quad}, \quad) = F$ as {\it false} (not {\it true}), the value $(\quad, \quad) = N$ as neither {\it true} nor {\it false}, and the value $(\cross{\quad}, \cross{\quad}) = B$ as both {\it true} and {\it false}. \\

 WF and $FOUR$ were  developed at roughly the same time in the 1970s, but attention has not previously been given to the relationship between the two systems. 

\section{The BF Calculus}\label{sec300}

The original motivation for BF was to have an operator $\cross{\quad}_i$ of order four rather than order two as in Spencer-Brown, so that $\cross{\cross{\quad}_i}_i = \cross{\quad}$, where $\cross{\quad}$ is the Spencer-Brown\textquotesingle s mark. Then  $\cross{\quad}_i$ can be considered a ``square root," analogous to the square root of negation and the square root of negative one. In so doing we have four basic values:
$\cross{\quad}_i$\ , $\cross{\cross{\quad}_i}_i = \cross{\quad}$\ , $\cross{\cross{\quad}}_i = \cross{\cross{\cross{\quad}_i}_i}_i $\ , and $\cross{\cross{\cross{\quad}}_i}_i =\  \cross{\cross{\quad}} = \quad \quad \quad  $.\\

We define the Brown-Four Calculus (BF) to share the same set of four of values  as the three calculi in the prior section. \\

\begin{definition}\label{7450} The four \emph{simple values} in BF are\\ (\quad, \quad), (\cross{\quad}, \quad), (\cross{\quad}, \cross{\quad}), and (\quad, \cross{\quad}).\\  
\end{definition}

BF shares the pairwise juxtaposition operation $A\ B$ as defined in Definition \ref{7300}. What distinguishes BF is the nesting operation $\cross{\quad}_i$\ , which we call the \emph{square root of negation}. \\

\begin{definition}\label{bf1000} Let $X = (a, b)$. We define the \emph{square root of negation} to be the nesting operation $\cross{\quad}_i$:\\ 

$\cross{X}_i = \cross{(a,\ b)}_i = (\cross{b}, a)$. \quad \quad \quad \quad \quad \quad \quad (SQRT)\\
\end{definition}

\noindent
As in \emph{Laws of Form}, we  identify the empty value with its representative, the empty parentheses:\\[-6pt]

%$$(\quad, \quad) = \quad \quad .$$ 

\begin{description}
\item $(\quad, \quad) = \quad \quad.$\\[-6pt]
\end{description}

\begin{figure}[htbp]
\centerline{ \includegraphics[width=2.6 in]{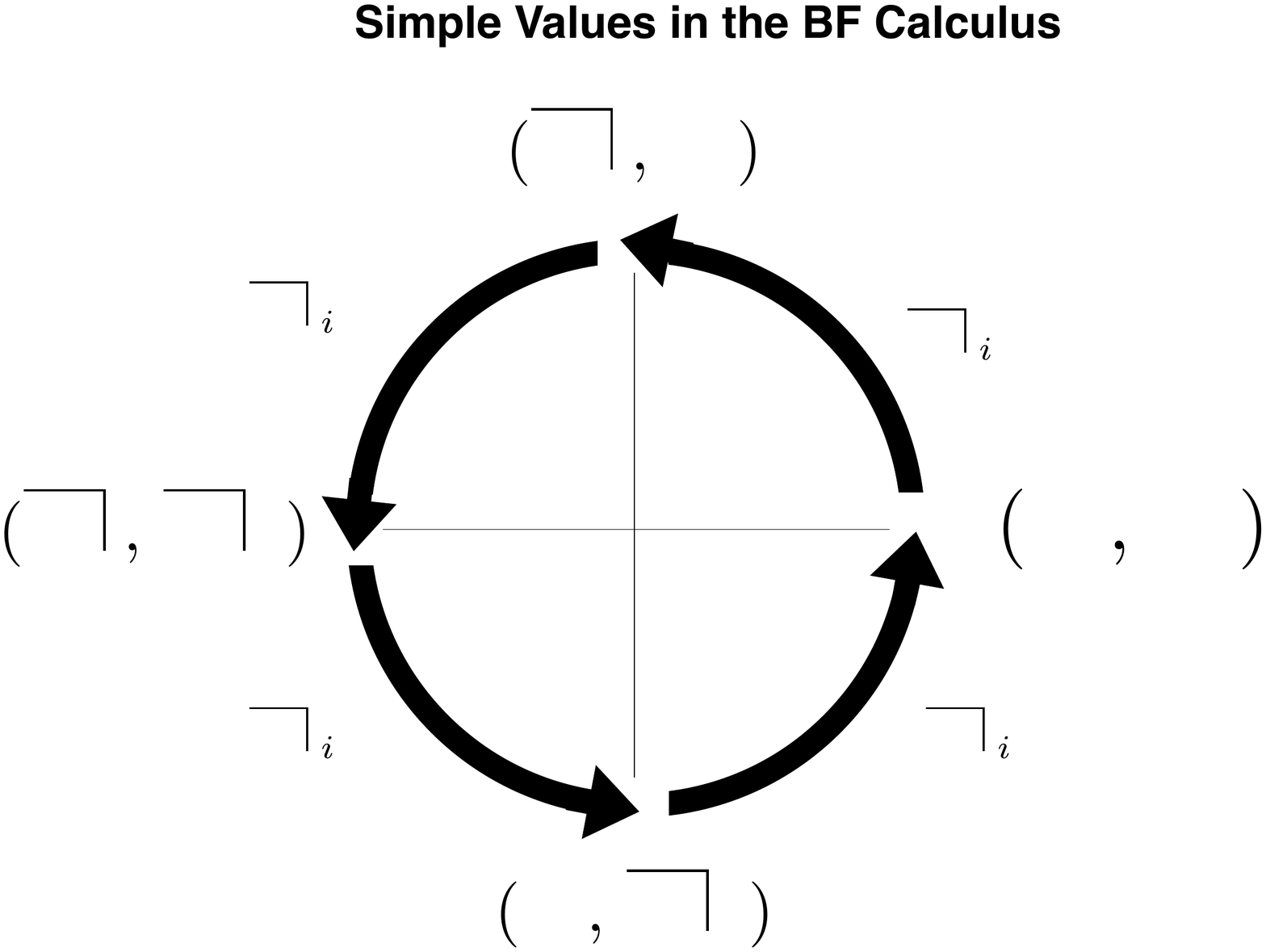}}
\caption{The Operation SQRT forms a 4-cycle.}
\label{fig2}
\end{figure}

\noindent
As illustrated in Figure \ref{fig2}, the four simple values are\\[-6pt] 
\begin{description}
\item $\cross{(\quad, \quad)}_i = (\cross{\quad}, \quad)$\\[-6pt]
\item $\cross{\cross{(\quad, \quad)}_i}_i = (\cross{\quad}, \cross{\quad})$\\[-6pt]
\item $\cross{\cross{\cross{(\quad, \quad)}_i}_i}_i = (\quad, \cross{\quad})$\\[-6pt]
\item $\cross{\cross{\cross{\cross{(\quad, \quad)}_i}_i}_i}_i = (\quad, \quad)$\\[-6pt]
\end{description}

\begin{figure}[htbp]
\centerline{ \includegraphics[width=2.0 in]{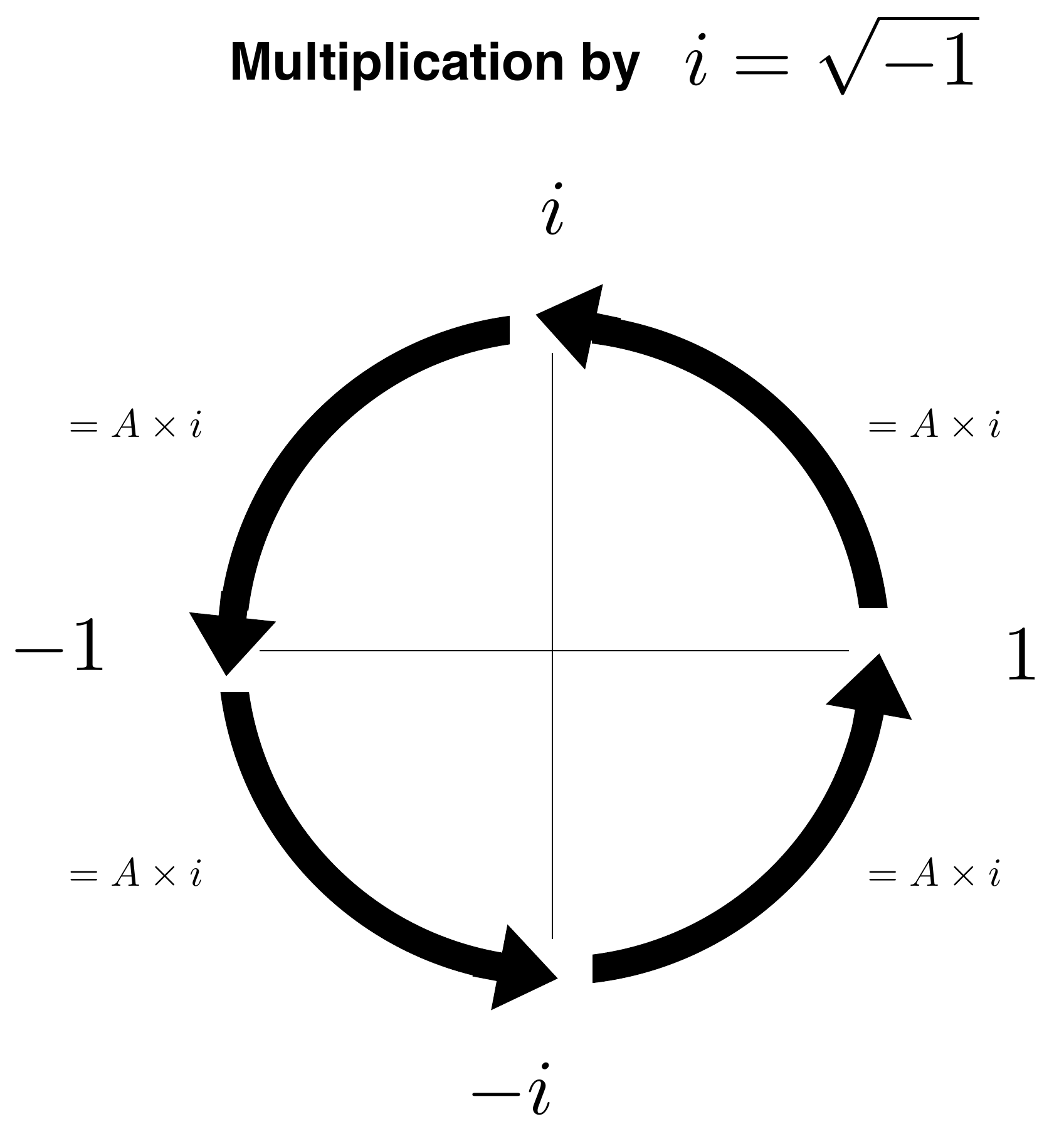}}
\caption{The operation SQRT is exactly analogous to multiplying by $i =\sqrt{-1}$}
\label{fig30}
\end{figure}

The operation $\cross{X}_i$ is analogous to multiplying unit values in the complex plane by $i = \sqrt{-1}$, as shown in Figure \ref{fig30}.  

The analogy is exact, based on the mapping $ (\quad, \quad) \rightarrow 1$, $ (\cross{\quad}, \cross{\quad}) \rightarrow -1$, $ (\cross{\quad}, \quad)\rightarrow i$, and $  (\quad, \cross{\quad}) \rightarrow -i$, and justifies our use of the term square root of negation.  \\

We make use of superscripts to indicate multiple nestings: $\cross{X}_i^2 = \cross{\cross{X}_i}_i $\ , and in general $\cross{X}_i^n = \cross{\cross{\cross{X}_i}_i \ \dots \  }_i$\ . In this notation the four simple values are given by applying operators to the unmarked state:
  \quad \quad, $\cross{\quad}_i$\ , $\cross{\quad}_i^2$\ , and $\cross{\quad}_i^3$.\\

\subsection{Calculating in BF}

An \emph{arithmetic} expression contains no variables. In BF each arithmetic expression simplifies to a simple expression. We will prove below that in BF each arithmetic expression simplifies to a simple expression. Note that any concatenation of operators applied to some ordered pair of expressions reduces to a new order pair of expressions each in the primary arithmetic of Spencer-Brown. Then each of these reduces to either a mark or a void. Thus every expression in BF does reduce to a simple value. We will below prove the uniqueness of this reduction.
$$(\cross{\quad}, \quad)$$\\[-24pt] $$(\cross{\cross{\cross{\quad}}\cross{\quad}}, \cross{\quad})$$\\[-24pt]
$$\cross{\cross{(\cross{\quad}, \cross{\quad})}_i\ \cross{\cross{\cross{(\cross{\quad}, \quad)}_i}_i}_i}_i$$\\[-24pt]

The first is already a simple expression. The second simplifies via two successive applications of Rule A2: $$(\cross{\cross{\cross{\quad}}\cross{\quad}}, \cross{\quad})   =(\cross{\cross{\quad}}, \cross{\quad})= (\quad, \cross{\quad}). $$
Simplification of the third proceeds in two steps. The first step uses only the operations SQRT (Def. \ref{bf1000}) and JUXT (Def. \ref{7300}): 
\begin{equation}\label{juxt_sqrt}
\begin{aligned}
&\cross{\cross{(\cross{\quad}, \cross{\quad})}_i\ \cross{\cross{\cross{(\cross{\quad}, \quad)}_i}_i}_i}_i \\[6pt]	
&= \cross{(\cross{\cross{\quad}}, \cross{\quad})\ \cross{\cross{(\cross{\quad}, \cross{\quad})}_i}_i}_i &[SQRT]\\[6pt]	
&= \cross{(\cross{\cross{\quad}}, \cross{\quad})\ (\cross{\cross{\quad}}, \cross{\cross{\quad}})}_i &[SQRT, twice]\\[6pt]	
&= \cross{(\cross{\cross{\quad}}\ \cross{\cross{\quad}}, \cross{\quad}\ \cross{\cross{\quad}})}_i &[JUXT]\\[6pt]	
&= (\cross{\cross{\quad}\ \cross{\cross{\quad}}}, \cross{\cross{\quad}}\ \cross{\cross{\quad}}) &[SQRT]
\end{aligned}
\end{equation}
The second step uses only rules A1 or A2 for the Primary Arithmetic: 
\begin{equation}\label{pa_simple}
\begin{aligned}
&(\cross{\cross{\quad}\ \cross{\cross{\quad}}}, \cross{\cross{\quad}}\ \cross{\cross{\quad}}) \\[6pt]
&= (\cross{\cross{\quad} }, \cross{\cross{\quad}}) &[A2, twice]\\[6pt]
&= (\quad,\quad) &[A2, twice]
\end{aligned}
\end{equation}

We will show that simplification of expressions is unique.\\ 

\begin{definition}\label{def_ums}
Let $E$ be an arithmetic expression. A \emph{Nested Marking Scheme} is a labeling of the spaces within $E$ with markings $M_0, M_1,\ M_2,\ M_3$ as follows. Begin at the left of E, labeling any constant values $(a,b)M_x: (\quad, \quad)M_0$, $(\cross{\quad}, \quad)M_1$, $(\cross{\quad}, \cross{\quad})M_2$, $(\quad, \cross{\quad})M_3$. Then label from left to right according to the following rules of nested marking:\\

\begin{enumerate}
    \item $\cross{M_x}_i  \rightarrow \cross{M_x}_i\ M_{x+1 (mod\ 4)}$
    \item $M_2\ M_x  \rightarrow M_2$, and $M_1 \ M_3 \rightarrow M_2$
    \item $M_x\ M_x \rightarrow M_x$
    \item $M_0\ M_x \rightarrow M_x.$
\end{enumerate}
Note that if the arithmetical expression has all its simple values written in terms of nested square roots of negation, then one need only start by labelling all unmarked spaces with $M_0.$ The simplified value of $E$ corresponds to the rightmost mark. Both examples below simplify to the value $(\cross{\quad},\quad) \longleftrightarrow M_1$. In actual practice we shall use the numerals $0,1,2,3$ for $M_0,M_1,M_2,M_3.$
\\[-6pt]
\end{definition}

\begin{example}
$$\cross{\cross{\cross{\cross{\cross{\quad}_i}_i}_i}_i}_i \longrightarrow \cross{\cross{\cross{\cross{\cross{0}_i\ 1}_i\ 2}_i\ 3}_i 0}_i\ 1$$
\end{example}

\begin{example}
$$\cross{\cross{\cross{\cross{(\cross{\ },\cross{\  } ) }_i (\cross{\  }, \  )  }_i}_i}_i  \rightarrow \cross{\cross{\cross{\cross{(\cross{\  },\cross{\  } )2 }_i 3\ (\cross{\  }, \  )1\ 2 }_i 3}_i0}_i1$$
\end{example}

\begin{theorem}\label{th315}
Every arithmetic expression $E$ in BF simplifies uniquely to one of the four simple values. 

\emph{Proof (Sketch)}: Per Theorem \ref{th100}, every arithmetic expression $(a,b)$ in $E$  simplifies uniquely to a simple value. Mark the resulting expression in accord with the Nested Marking Scheme per Definition \ref{def_ums}. Then we claim the labeled output (far right) corresponds to the simplified value of $E$.  We observe the operations SQRT and JUXT each leave the marking scheme unchanged. But any simplification of $E$ must occur via a sequence of steps JUXT and SQRT. Consequently, $E$ cannot simplify to more than one simple value. 
\end {theorem}

\subsection{Algebraic Expressions}

\begin{definition}\label{def_var} As in Definition~\ref{exp_d1} an algebraic expression in BF is an expression that may have variables standing for possible arithmetic values. There will be, in such an expression, some number   ($\ge 0$) of variables $X = (x_1, x_2)$, where $x_1$ and $x_2$ are variables in the Primary Algebra. Note that an algebraic expression will be reducible to a pair of primary algebraic expressions.\\
\end{definition}

\begin{definition}\label{74400} \emph{Equality}. Two algebraic expressions $F$ and $G$ are \emph{arithmetically} equivalent if they are equal in every case of choosing one of the four simple values for each variable. They are \emph{demonstrably} equivalent if $F = G$ can be demonstrated as a series of steps applying JUXT, SQRT, and the axioms of the Primary Algebra. Note this means one can apply the rules of the primary algebra independently to two algebraic expressions $a$ and $b$ that appear in a pair $(a,b).$ 
\\[-6pt]
\end{definition}

We give several important algebraic demonstrations below. \\[-6pt]

 \noindent
$B1.  \ \ \cross{\quad} \ A = \cross{\quad}$\quad \quad \quad \ (Integration) \\[-6pt]

\noindent
$\ \ \cross{\quad} \ A$\\[5pt]
\begin{tabular}{@{\hspace{3ex}}l l }
= $(\cross{\quad}, \cross{\quad})\ (a_1, a_2)$ &[Notation]\\[2pt]
= $(\cross{\quad}\ a_1, \cross{\quad}\ a_2)$ &[JUXT (Def. \ref{7300})]\\[2pt]
= $(\cross{\quad}, \cross{\quad})$ &[B2, Integration in the PA]\\[2pt]
= $\cross{\quad}$ &[Notation]\\[6pt]
\end{tabular}

\noindent 
$B2. \ \ \   \cross{A}_i^4 = A$\quad \quad \quad \quad \quad (Reflexion)\\

$\ \ \cross{(a_1, a_2)}_i^4$\\[5pt]
\begin{tabular}{@{\hspace{3ex}}l l }
= $\cross{(\cross{a_2}, a_1)}_i^3$ &[SQRT (Def. \ref{bf1000})]\\[2pt]
= $\cross{(\cross{a_1}, \cross{a_2})}_i^2$ &[SQRT]\\[2pt]

= $\cross{(\cross{\cross{a_2}}, \cross{a_1})}_i$ &[SQRT]\\[2pt]
= $(\cross{\cross{a_1}}, \cross{\cross{a_2}})$ &[SQRT]\\[2pt]
= $(a_1, a_2)$ &[B2, Reflexion in the PA]\\[2pt]
= $A$ &[Notation]\\[6pt]
\end{tabular}

 \noindent 
 B3\#. $\ \ \cross{\cross{A}_i\ B}_i\ C =   \cross{\cross{A\ C}_i\ B}_i\ C$\ \  (Split Generation)\\[6 pt]
$\cross{\cross{A}_i\ B}_i\ C$\\[5pt]
\begin{tabular}{@{\hspace{4ex}}l l }
= $\cross{\cross{(a_1, a_2)}_i\ (b_1, b_2)}_i\ (c_1, c_2)\ $ &[Notation]\\[3pt]
= $\cross{(\cross{a_2}, a_1)\ (b_1, b_2)}_i\ (c_1, c_2)\ $ &[SQRT]\\[3pt]
= $\cross{(\cross{a_2}\ b_1,\ a_1\ b_2)}_i\ (c_1, c_2)\ $ &[JUXT]\\[3pt]
= ${(\cross{a_1\ b_2},\ \cross{a_2}\ b_1)}\ \ (c_1, c_2)\ $ &[SQRT]\\[3pt]

= $(\cross{a_1\ b_2} c_1,\ \cross{a_2}\ b_1  c_2)\ $ &[JUXT]\\[3pt]

= $(\cross{a_1\ c_1\ b_2}\ c_1,\ \cross{a_2 \ c_2}\ b_1\ c_2)\ $ &[Gen. in PA]\\[3pt]
= $(\cross{a_1\ c_1\ b_2},\ \cross{a_2 \ c_2}\ b_1)\ (c_1, c_2)$ &[JUXT]\\[3pt]
= $\cross{(\cross{a_2 \ c_2}\ b_1,\ a_1\ c_1\ b_2)}_i\ (c_1, c_2)$ &[SQRT]\\[3pt]
%\end{tabular}
%\begin{tabular}{@{\hspace{4ex}}l l }
= $\cross{(\cross{a_2 \ c_2},\ a_1\ c_1)\ (b _1,\ b_2)}_i\ (c_1, c_2)$ &[JUXT]\\[3pt]
= $\cross{\cross{(a_1\ c_1,\ a_2 \ c_2)}_i\ (b_1,\ b_2)}_i\ (c_1, c_2)$ &[SQRT]\\[3pt]
\end{tabular}
\begin{tabular}{@{\hspace{4ex}}l l }
= $\cross{\cross{(a_1, a_2)\ (c_1, c_2)}_i\ (b_1,\ b_2)}_i\ (c_1, c_2)$ &[JUXT]\\[3pt]
= $\cross{\cross{A\ C}_i\ B}_i\ C$ &[Notation]\\[12pt]
\end{tabular}

\noindent B11. $\cross{\cross{A}_i^3\ \cross{B}_i^3}_i\ C = \cross{\cross{A\ C}_i^3\ \cross{B\ C}_i^3}_i\ $\\[-6pt]

\noindent B12. $\cross{\cross{A}_i\ \cross{B}_i}_i^3\ C = \cross{\cross{A\ C}_i\ \cross{B\ C}_i}_i^3\ $.\\[-6pt]

We leave the final demonstrations for the reader. B3\# is a new result that generalizes the form of Generation (Axiom B3). Consequences B4-B10 (Primary Algebra) are each valid in BF. B11 and B12 are new distribution laws, closely related to the meet and join operations for bilattices in Section \ref{sec_v}.  
\subsection{Axiomatic Completeness}

\begin{theorem}\label{th4050}  Let  $F$ be an algebraic expression in BF. Then $F$ is demonstrably equivalent to an expression $\hat{F} = (F_1, F_2)$ where $F_1$, and $F_2$ are expressions in the  primary algebra.\\[-6pt]

\emph{Proof}: 1. If $F$ has depth 0 then , $F = (a_1, a_2)(a_3,a_4)\\ \dots (a_{2k-1},a_{2k})$ $= (a_1a_3\dots a_{2k-1}, a_2a_4\dots a_{2k}) = (F_1, F_2)$, where $F_1$ and $F_2$ are both expressions in the PA. \\[-6pt]

2. Otherwise assume by induction the theorem to be true for expressions of depth $\le n-1$, and that $F$ has depth $n$. Then $F$ can be reduced to depth $n-1$  by applying the result from 1 to each of its deepest spaces, followed by applying the SQRT operation to each such deepest space. Being reduced to depth $n-1$, it then follows that $F = \hat{F} = (F_1, F_2)$. \\
\end{theorem}

\begin{theorem}\label{th4051}
\emph{Axiomatic Completeness}. Let $F = G$ be an arithmetical equivalence. Then $F = G$ is also demonstrable. \\[-6pt]

\emph{Proof:}
  By Theorem \ref{th4050}, $F = \hat{F} = (F_1, F_2)$ and $G = \hat{G} = (G_1, G_2)$ are demonstrable equalities, so therefore, $(F_1, F_2) = (G_1, G_2)$. Consequently, it must also be the case that $F_1 = G_1$ and $F_2 = G_2$ are both arithmetical equalities, since otherwise by examination of truth tables (choices of simple values for the variables) we would necessarily find an exception to Theorem \ref{th315}. Then by Theorem \ref{th200}, $F_1 = G_1$ and $F_2 = G_2$ are also demonstrable as expressions in the Primary Algebra, and so therefore is $F = G$.
\end{theorem}

\begin{remark}\label{completeness} See \cite{bf7, bf8, bf81} for a proof that the three axioms B1, B2, and B3\# form a complete basis for the BF calculus, without making reference to the primary algebra.  \\[-6pt]
\end{remark}

\begin{description}
\item [B1.] $\cross{\cross{\cross{\cross{A}_i}_i}_i}_i = A$\quad \quad \quad \quad \quad \ \quad   (Reflexion)\\
\item [B2.] $A\ \cross{\cross{\quad}_i}_i= \cross{\cross{\quad}_i}_i$\quad \ \  \quad \quad  \quad (Integration)\\
\item [B3\#.] $\cross{\cross{A}_i\ B}_i\ C = \cross{\cross{AC}_i\ B}_i\ C$\quad (Split Generation)\\[-6pt]
\end{description}

\noindent
Note that $\cross{\ \ } = \cross{\cross{\ \ }_i}_i$.\\ 

\noindent
We can easily show that B3\# implies B3.\\[-6pt]

$ \cross{A} B$\\[5pt]
\begin{tabular}{@{\hspace{8ex}}l l }
= $ \cross{\cross{A}_i}_i\ B$ &\\[2pt]
= $ \cross{\cross{AB}_i}_i\ B$ &[B3\#]\\[2pt]
= $ \cross{AB}\ B$. \\[6pt]
\end{tabular}

\noindent
B1, B2, and B3\# implicate the rules for BF arithmetic. Previously we had \small (\cross{ \quad}, \quad)\ (\quad, \cross{\quad}) = (\cross{\quad}, \cross{\quad})\normalsize. Now we find as we have seen by pairs \\[-6pt]

$\cross{\quad}_i\ \cross{\quad}_i^3 $\\[5pt]
\begin{tabular}{@{\hspace{8ex}}l l }
= $\cross{\quad}_i\ \cross{\cross{\quad}_i}$&[By Definition]\\[2pt]
= $\cross{\quad}_i\ \cross{\quad}$&[B3, Generation]\\[2pt]
= $\cross{\quad}$&[B2, Integration].\\[2pt]
\end{tabular}

\section{BF and the Bilattice FOUR}\label{sec_v}

In this section, we interpret the BF calculus as a 4-valued \emph{bilattice}. The bilattice concept was initially developed by Nuel Belnap \cite{bf10}, and refined by  Ginsburg \cite{bf11} and Fitting \cite{bf12}. The simplest bilattice is known as Belnap's FOUR. A central idea of this paper is showing that the bilattice can expressed using only two (rather than six) primitive operations. 

\subsection{The Bilattice FOUR}

As originally conceived, $FOUR$ consists of two sets of meet and join operations ($\wedge, \vee, \otimes, \oplus$) and one negation operation ($\neg$), plus the four  values {\it true} ($T$), {\it false} ($F$), Neither {\it true}  nor {\it false} $(N)$, and Both {\it true} and {\it false} $(B)$. For this paper, we will consider a bilattice to include a second negation operation ($!$) called ``conflation," as introduced by Fitting. \cite{bf12}\\[-6pt]

The four values are often represented as ordered pairs, with the first entry representing falseness, and the second, truth. Thus, $F =(1,0)$ means {\it false} and not {\it true}, $T = (0,1)$ {\it true} and not {\it false}, $N =(0,0)$ neither, and $B =(1,1)$ both. We map these values to BF as follows: \\[-18pt]

\begin{equation}
\begin{aligned}
&F \longleftrightarrow (\cross{\quad}, \quad)	&T \longleftrightarrow (\  \quad, \cross{\quad})\\
&N \longleftrightarrow (\quad, \quad)	&\ B \longleftrightarrow (\cross{\quad}, \cross{\quad})
\end{aligned}
\end{equation}

\subsection{The Six Bilattice Operations}

\begin{figure}[htbp]
\centerline{ \includegraphics[width=1.8in]{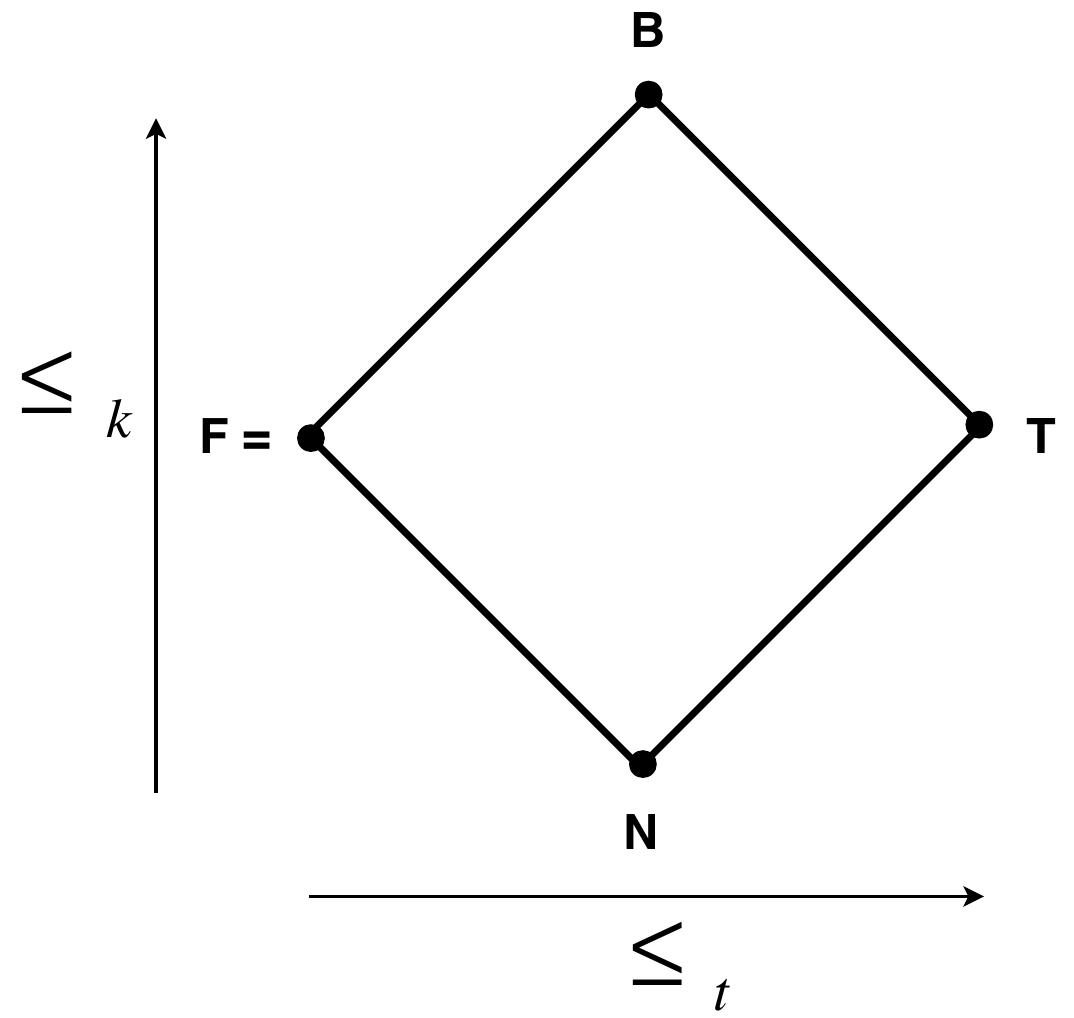}}
\caption{Lattice Diagram for the Bilattice FOUR}
\label{fig3}
\end{figure}

The meet and join operations are defined as the greatest lower (glb) and least upper (lub) bounds on the two partial orderings  $\le_t$ and $\le_k$ in the lattice diagram in Figure \ref{fig3}. The ``truth ordering" $\le_t$ defines $A \wedge B = glb_{t}(A,B)$ and $A \vee B = lub_{t}(A,B)$, and in the ``knowledge ordering" $\le_k$ defines $A \otimes B = glb_{k}(A,B)$ and $A \oplus B = lub_{k}(A,B).$\\

Negation operations for the two orderings are defined \\[-15pt]

\begin{equation}
\begin{aligned}
&\neg A = \neg(a,b) = (b,a)	&(Ordering \le_t)\\
&!A =\ !(a,b) = (\sim b,\sim a)	&(Ordering \le_k).\\
\end{aligned}
\end{equation}
\noindent 
Output values for the first are $\neg T=F$, $\neg F=T$, $\neg N=N$, and $\neg B=B$; and for the second, $!T=T$, $!F=F$, $!N=B$, and $!B=N$. The two sets of operations form separate, interrelated DeMorgan Algebras. Fitting observes that $\neg $ and $!$ mutually commute, with $\neg !A =\ !\neg A =\ \sim A$. \cite{bf12}

\subsection{Interpretation}\label{interpret}

We interpret the bilattice operations in BF as follows: 
\begin{eqnarray*}
 &Bilattice	 &BF\\
  A \vee B &= lub_{t}(A,B)	 &\cross{\cross{A}_i^3 \  \cross{B}_i^3}_i\\
 A \wedge B &= glb_{t}(A,B)	& \cross{\cross{A}_i \  \cross{B}_i}_i^3\\
 \neg A &= (a_2, a_1) &\cross{\cross{A}_i\ \cross{\quad}_i^3\ } \ \cross{\cross{A}_i^3\ \cross{\quad}_i\ }
 \end{eqnarray*}
 The three operations above define a DeMorgan Algebra on the t-ordering in the lattice diagram in Figure \ref{fig3} that corresponds to Belnaps FOUR, while the three operations below form a DeMorgan Algebra on the k-ordering that corresponds with the Waveform Algebra WF. 
 \begin{eqnarray*}
 A \oplus B &= lub_{k}(A,B)	& A\ B\\
 A \otimes B &= glb_{k}(A,B)	& \cross{\cross{A}\ \cross{B}}\\
 !A &= ( \sim a_2, \sim a_1)& \cross{\cross{A}_i\ \cross{\quad}_i\ } \ \cross{\cross{A}_i^3\ \cross{\quad}_i^3\ } 
\end{eqnarray*}

\noindent  
Inspecting the following truth tables is useful in understanding the definitions of these expressions. \\ 

\begin{tabular}{l c}
\begin{tabular}{|c|c|c|c|c|}
\hline 
\multicolumn{5}{|l|}{$\mathbf{Table\ 1:}$ \scriptsize  $\cross{\cross{A}_i^3 \  \cross{B}_i^3}_i$  } \\
\hline
$A \vee B$		&$\boldsymbol{1}$    &$\boldsymbol{0}$  &$\boldsymbol{2}$ &$\boldsymbol{3}$ \\
\hline
 $\boldsymbol{1}$     	&$1$   &$0$ &$2$ &$3$ \\
\hline
$\boldsymbol{0}$    &$0$ &$0$ &$3$ &$3$\\
\hline
$\boldsymbol{2}$	&$2$  &$3$ &2 &$3$  \\
\hline
$\boldsymbol{3}$	&$3$  &$3$ &$3$ &$3$  \\
\hline
\end{tabular}\ \ 

\begin{tabular}{|c|c|c|c|c|}
\hline 
\multicolumn{5}{|l|}{$\mathbf{Table\ 3:}$ \large  $A \ B$  } \\
\hline
$A \oplus B$		&$\boldsymbol{0}$    &$\boldsymbol{1}$  &$\boldsymbol{3}$ &$\boldsymbol{2}$ \\
\hline
 $\boldsymbol{0}$     	&$0$   &$1$ &$3$ &$2$ \\
\hline
$\boldsymbol{1}$    &$1$ &$1$ &$2$ &$2$\\
\hline
$\boldsymbol{3}$	&$3$  &$2$ &3 &$2$  \\
\hline
$\boldsymbol{2}$	&$2$  &$2$ &$2$ &$2$  \\
\hline
\end{tabular}

\end{tabular}\\\\

\begin{tabular}{l c}
\begin{tabular}{|c|c|c|c|c|}
\hline 
\multicolumn{5}{|l|}{$\mathbf{Table\ 2:}$ \scriptsize  $\cross{\cross{A}_i \  \cross{B}_i}_i^3$  } \\
\hline
$A \wedge B$		&$\boldsymbol{1}$    &$\boldsymbol{0}$  &$\boldsymbol{2}$ &$\boldsymbol{3}$ \\
\hline
 $\boldsymbol{1}$     	&$1$   &$1$ &$1$ &$1$ \\
\hline
$\boldsymbol{0}$    &$1$ &$0$ &$1$ &$0$\\
\hline
$\boldsymbol{2}$	&$1$  &$1$ &2 &$2$  \\
\hline
$\boldsymbol{3}$	&$1$  &$0$ &$2$ &$3$  \\
\hline
\end{tabular}\ \ 

\begin{tabular}{|c|c|c|c|c|}
\hline 
\multicolumn{5}{|l|}{$\mathbf{Table\ 4:}$  \scriptsize $\cross{\cross{A} \ \cross{B}}$  } \\
\hline
$A \otimes B$		&$\boldsymbol{0}$    &$\boldsymbol{1}$  &$\boldsymbol{3}$ &$\boldsymbol{2}$ \\
\hline
 $\boldsymbol{0}$     	&$0$   &$0$ &$0$ &$0$ \\
\hline
$\boldsymbol{1}$    &$0$ &$1$ &$0$ &$1$\\
\hline
$\boldsymbol{3}$	&$0$  &$0$ &3 &$3$  \\
\hline
$\boldsymbol{2}$	&$0$  &$1$ &$3$ &$2$  \\
\hline
\end{tabular}

\end{tabular}\\

\noindent
The four values are represented in the tables as follows: \\[-6pt]
\begin{description}\small
\item $0 = (\quad, \quad)$\\[-6pt]
\item $1 = \cross{(\quad, \quad)}_i = (\cross{\quad}, \quad)$\\[-6pt]
\item $2 = \cross{(\quad, \quad)}_i^2 = (\cross{\quad}, \cross{\quad})$\\[-6pt]
\item $3 = \cross{(\quad, \quad)}_i^3 = (\quad, \cross{\quad})$\\
\end{description} 

\noindent
Tables 1 and 2 are based on the ordering in the t-lattice, while tables 3 and 4 are based on the ordering in the k-lattice. Values in the table are arranged to emphasize the isomorphic relation between tables 1/3 and 2/4.\\ 

The following demonstrations clarify the basis for our definitions of the DeMorgan negation operations $\neg$ and $!$.\\

\small
\noindent
\emph{T-Ordering}: $\neg A = (a_2, a_1)$\\[6pt]
$\cross{\cross{A}_i\ \cross{\quad}_i^3\ } \ \cross{\cross{A}_i^3\ \cross{\quad}_i\ }$\\[5pt]
\begin{tabular}{@{\hspace{3ex}}l l }
= $\cross{\cross{(a_1, a_2)}_i\ (\quad, \cross{\quad}) } \ \cross{\cross{(a_1, a_2)}_i^3\ (\cross{\quad}, \quad) }$ &[Notation]\\[3pt]
= $\cross{(\cross{a_2}, a_1)\ (\quad, \cross{\quad}) } \ \cross{(a_2, \cross{a_1})\ (\cross{\quad}, \quad) }$ &[SQRT]\\[3pt]
= $\cross{(\cross{a_2}, \cross{\quad}) } \ \cross{(\cross{\quad}, \cross{a_1}) }$ &[JUXT]\\[3pt]
= $(a_2, \quad)  \ \ (\quad, a_1) $ &[SQRT 2X]\\[3pt]
= $(a_2, a_1)$ &[JUXT]\\[3pt]\\
\end{tabular}
\normalsize

\noindent
\small 
\emph{K-Ordering}: $!A = (\cross{a_2}, \cross{a_1})$\\[6pt]
$\cross{\cross{A}_i\ \cross{\quad}_i\ } \ \cross{\cross{A}_i^3\ \cross{\quad}_i^3\ }$\\[5pt]
\begin{tabular}{@{\hspace{3ex}}l l }
= $\cross{\cross{\cross{A}_i\ \cross{\quad}_i^3\ } \ \cross{\cross{A}_i^3\ \cross{\quad}_i\ }}$ &[B10, CrossTransposition]\\[3pt]
= $\cross{(a_2, a_1)} $ &[By the Prior Demonstration]\\[3pt]
= $(\cross{a_2}, \cross{a_1})$ &[SQRT 2X]\\[6pt]
\end{tabular}
\normalsize

\noindent
From the above we see that $!\neg A = \neg !A =\ \sim A = \cross{A}$, as originally observed by Fitting, who regards a bilattice to be the ``lattice product" of the two DeMorgan algebras. With  inclusion of the identity operation $A$, the three negation operations $\neg A,\ \sim A$, and $!A$  form a group (under composition of expressions) that is isomorphic to Klein Four. \\

There are also additional symmetries. Consider the following expressions. \\

\small
\noindent
$\cross{\cross{A}\ \cross{\quad}_i^3\ } \ \cross{A \ \cross{\quad}_i\ }$\quad \quad \quad \quad \normalsize (Mark Right)\small\\[5pt]
\begin{tabular}{@{\hspace{3ex}}l l }
= $\cross{\cross{(a_1, a_2)}\ (\quad, \cross{\quad}) } \ \cross{(a_1, a_2)\ (\cross{\quad}, \quad) }$ &[Notation]\\[6pt]
= $\cross{(\cross{a_1}, \cross{a_2})\ (\quad, \cross{\quad}) } \ \cross{(a_1, a_2)\ (\cross{\quad}, \quad) }$ &[SQRT 2X]\\[6pt]
= $\cross{(\cross{a_1}, \cross{\quad}) } \ \cross{(\cross{\quad}, a_2) }$ &[JUXT]\\[6pt]
= $(a_1, \quad)  \ \ (\quad, \cross{a_2}) $ &[SQRT 2X]\\[6pt]
= $(a_1, \cross{a_2})$ &[JUXT]\\[3pt]\\
\end{tabular}\\
\normalsize

\small
\noindent
$\cross{A\ \cross{\quad}_i^3\ } \  \cross{\cross{A}\ \cross{\quad}_i\ }$\quad \quad \quad \quad \normalsize (Mark Left)\small\\[5pt]
\begin{tabular}{@{\hspace{3ex}}l l }
= $\cross{(a_1, a_2)\ (\quad, \cross{\quad}) } \ \cross{\cross{(a_1, a_2)}\ (\cross{\quad}, \quad) }$ &[Notation]\\[6pt]
= $\cross{(a_1,a_2)\ (\quad, \cross{\quad}) } \ \cross{(\cross{a_1}, \cross{a_2})\ (\cross{\quad}, \quad) }$ &[SQRT 2X]\\[6pt]
= $\cross{a_1, \cross{\quad}) } \ \cross{(\cross{\quad}, \cross{a_2}) }$ &[JUXT]\\[6pt]
= $(\cross{a_1}, \quad)  \ \ (\quad, a_2) $ &[SQRT 2X]\\[6pt]
= $(\cross{a_1}, a_2)$ &[JUXT]\\[3pt]\\
\end{tabular}
\normalsize

Together with the expression $\cross{A}$ and the identity expression $A$, Mark Right and Mark Left also form a group that is isomorphic to the Klein Four Group. Furthermore, the following eight expressions also form a group:\\[-6pt] 

\begin{enumerate}
\item The expressions $A$, $\cross{A}_i$, $\cross{A}$, and $\cross{A}_i^3$.\\ 
\item The two DeMorgan Negations\\
\begin{description}
\item $\cross{\cross{A}_i\ \cross{\quad}_i\ } \ \cross{\cross{A}_i^3\ \cross{\quad}_i^3\ }$\\
\item $\cross{\cross{A}_i^3\ \cross{\quad}_i\ } \ \cross{\cross{A}_i\ \cross{\quad}_i^3\ }$\\
\end{description}
\item The expressions Mark Right and Mark Left\\
\begin{description}
\item $\cross{\cross{A}\ \cross{\quad}_i^3\ } \ \cross{A \ \cross{\quad}_i\ }$\\
\item $\cross{A\ \cross{\quad}_i^3\ } \  \cross{\cross{A}\ \cross{\quad}_i\ }$.\\
\end{description}
\end{enumerate}

This group is formed by composition of expressions and is isomorphic to the Dihedral Group of order 8, also known as the Symmetries of the Square. The expressions $A$, $\cross{A}_i$, $\cross{A}$, and $\cross{A}_i^3$ correspond to rotations, while the two DeMorgan Negations and Left/Right Mark correspond to reflections. The strongly motivated reader is encouraged to verify this result.  

\subsection{The BF Represention of the Bilattice}

There is much to say about the characteristics revealed by the BF representation of the bilattice.  We have succeeded in deriving a fundamentally new representation of the bilattice -- in the form of the BF calculus. Remarkably, the number of primitive operations is reduced from six in the bilattice to merely two: the juxtaposition $X\ Y$, and imaginary containment,  $\cross{X}_i\ $ (as anticipated but developed differently in \cite{bf20, bf21}). BF is both axiomatically and functionally complete, and provides a very compact axiomatization for the bilattice. 

One promising area for extending research is the application of the BF formulation of bilattice operations to modal logics, following the work of Alexander Karpenko\cite{bf14}. Our preliminary results in this area show significant promise in expressing a range of modal logics using BF forms, including the interesting logics GL and GRZ associated with the ``Logic of Provability" of George Boolos \cite{bf15}.   

\subsection{Beyond the Bilattice}

\begin{figure}[htbp]
\centerline{ \includegraphics[width=2.0 in]{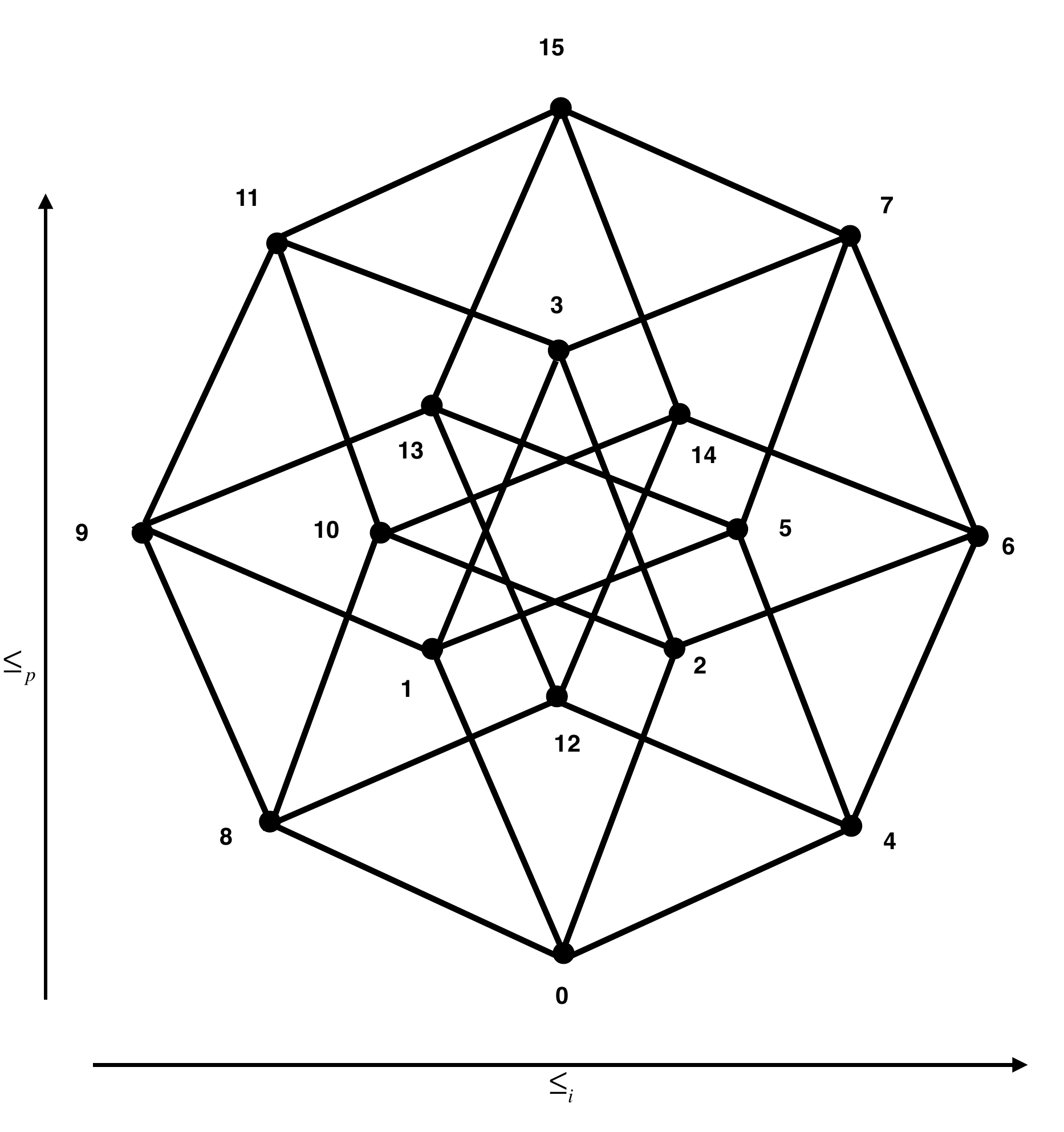}}
\caption{Lattice Diagram Defining Juxtaposition in the $Q$ Calculus}
\label{fig5}
\end{figure}

Based on additional work, we have shown that the concept of an imaginary operator extends to a 16-valued quaternion-based system $Q$ by adding the operator/values $\cross{\quad}_j$ and $\cross{\quad}_k$. In this quaternion system, $\cross{\cross{\quad}_i}_i = \cross{\cross{\quad}_j}_j = \cross{\cross{\quad}_k}_k = \cross{\cross{\cross{\quad}_i}_j}_k = \cross{\quad}$. The lattice of meets and joins that defines the juxtaposition operation forms a hypercube (Figure \ref{fig5}). 

We note that a similar construction of a 16-valued multi-lattice has previously been described by Shramko and Wansing\cite{bf16}. Our Calculus $Q$, based on the Quaternions, is conceptually simpler, and provides access to the axioms and consequences of both BF and the Primary Algebra. \\

We have also established that we can generalize the  calculus to the order $2^n$ for arbitrary $n$ using the rotational operation

    $$\cross{(x_1, x_2, \dots x_n)}_r = (\cross{x_n}, x_1, x_2, \dots , x_{n-1}). $$
  
The resulting system is functional complete and can be axiomatized by applying BF operations on adjacent indexed pairs. 

\subsection{Braiding}
\begin{figure}[htbp]
\centerline{ \includegraphics[width=2.0 in]{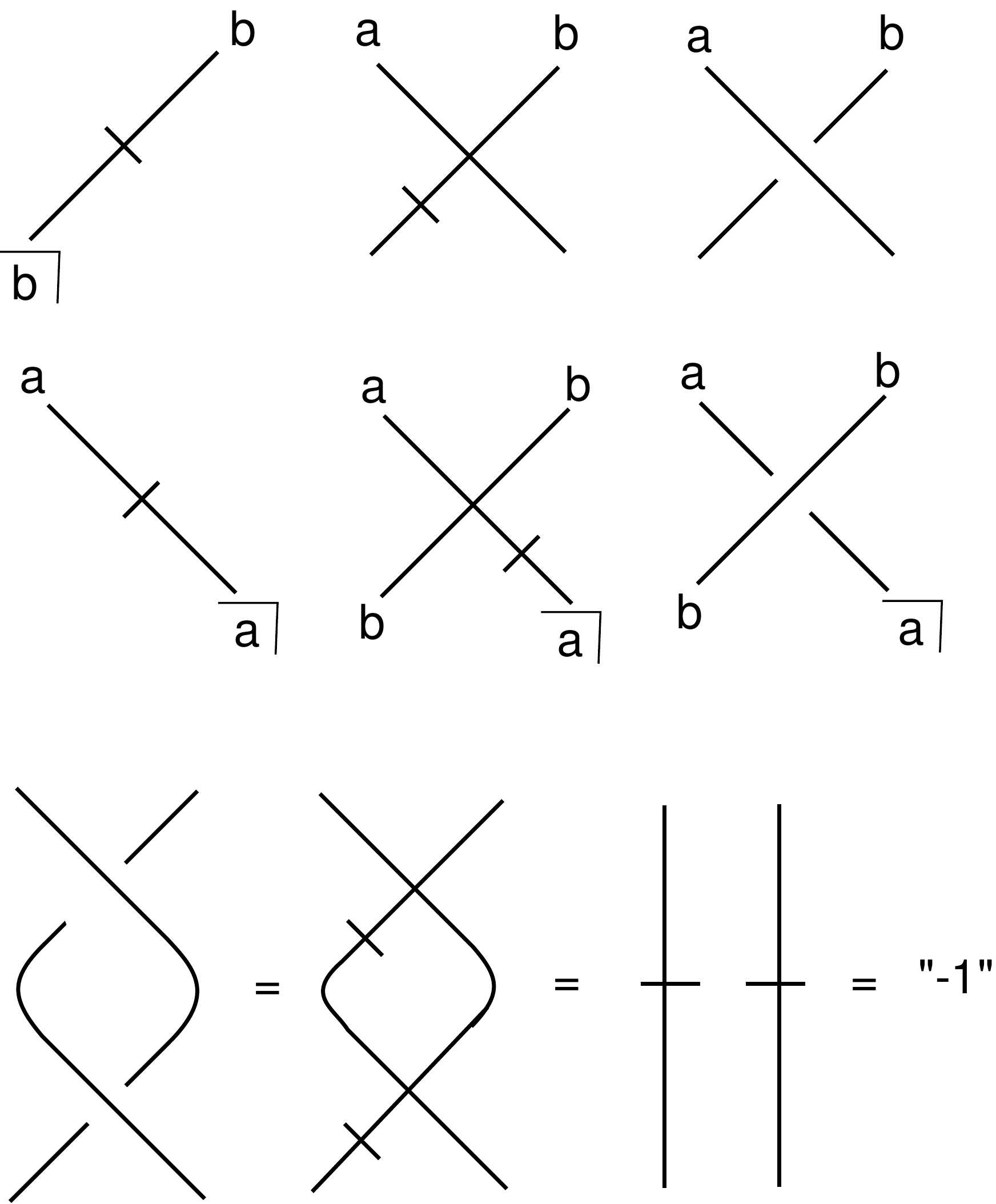}}
\caption{Basic Braiding Operations}
\label{basicbraid}
\end{figure}

\begin{figure}[htbp]
\centerline{ \includegraphics[width=2.0 in]{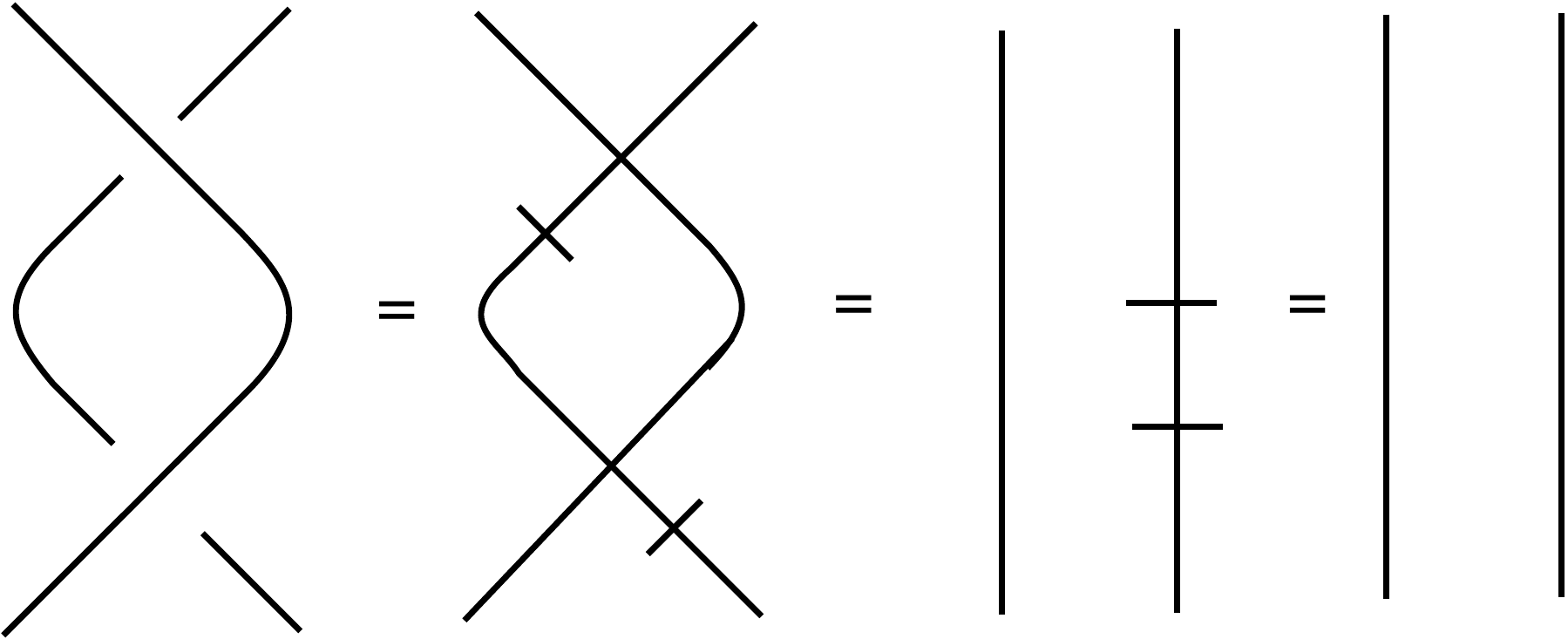}}
\caption{Braid Inverse}
\label{braidinverse}
\end{figure}

\begin{figure}[htbp]
\centerline{ \includegraphics[width=2.0 in]{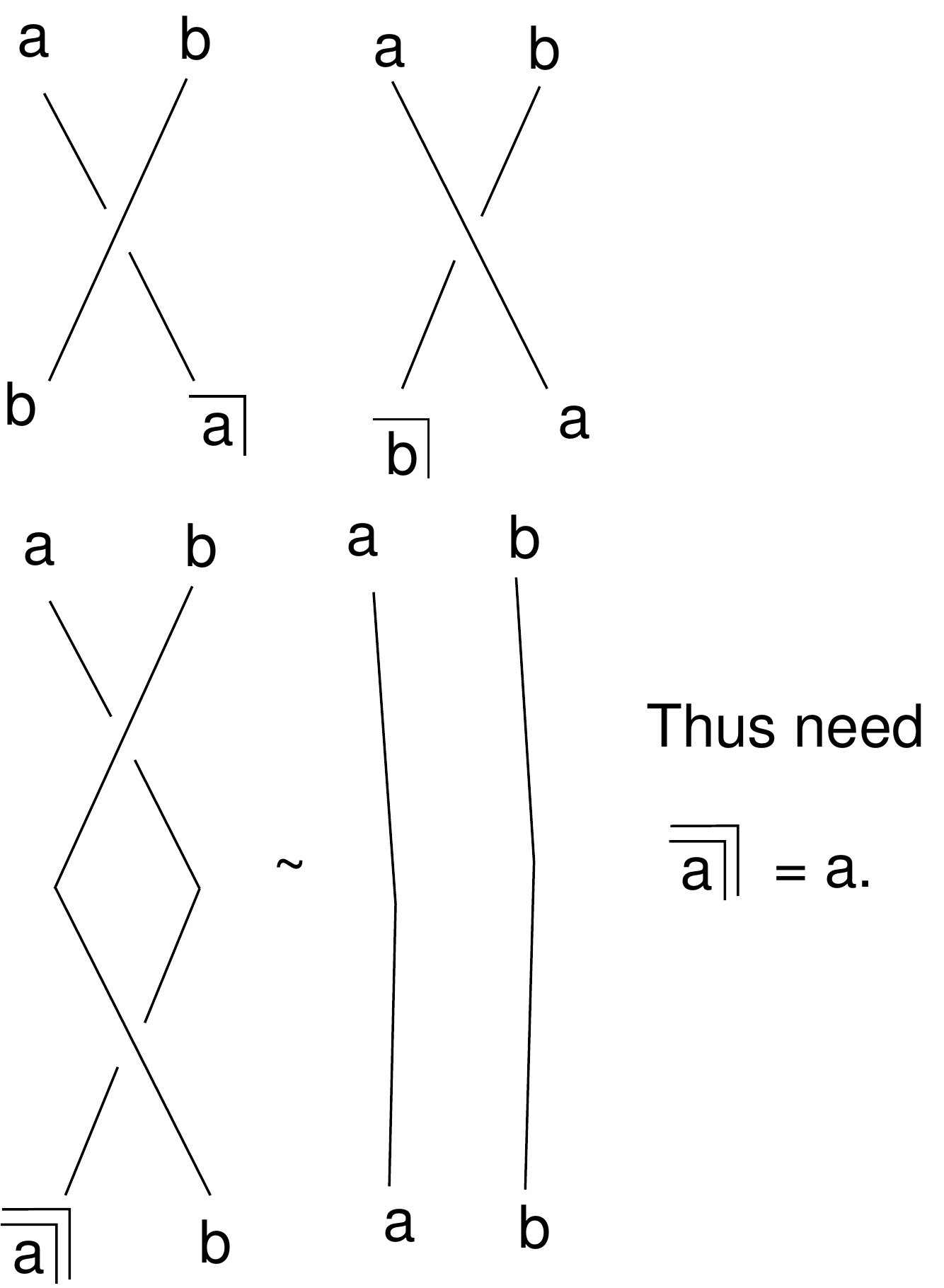}}
\caption{Representing imaginary marks as a braiding operations.}
\label{fig6}
\end{figure}

Our concept of these operations and of the square root of unity has a natural representation in terms of braiding operations.  In this representation, a pair (a,b) is represented as two strings, and the the imaginary operations $\cross{(a,b)}_i = (\cross{b}, a)$ and $\cross{(a,b)}_i^3 = (b, \cross{a})$ are represented so that the over-crossing string acts as the mark (operating on the label of the undercrossing string) and the crossing of the strings acts as the permutation. Note that in Figure  \ref{basicbraid}, Figure  \ref{braidinverse}  and Figure \ref{fig6} we have illustrated the nesting of the operators $\cross{\quad}_i$  and $\cross{\cross{\quad}_i}$ (twist and reverse twist) and we find that it will be the identity operator just so long as $\cross{\cross{a}} = a$ for any  $a$.\\

In these figures, we have indicated the braid itself as an operator that has the over-crossing acting as the mark, and we have also indicate the actual representation of the braid as a combination of the crossing operator and the permutation operator (a transposition) of the two strands. In Figure  \ref{basicbraid} we illustrate how the square of the two-strand braiding operator has the same structure as two parallel lines, each endowed with a crossing operator. This can be regarded as the ``-1" for the square of this operator and so we see that the fourth power of the operator is equal to the identity. See also  Figure  \ref{fig9} for another point of view on this fact.\\

\begin{figure}[htbp]
\centerline{ \includegraphics[width=2.0 in]{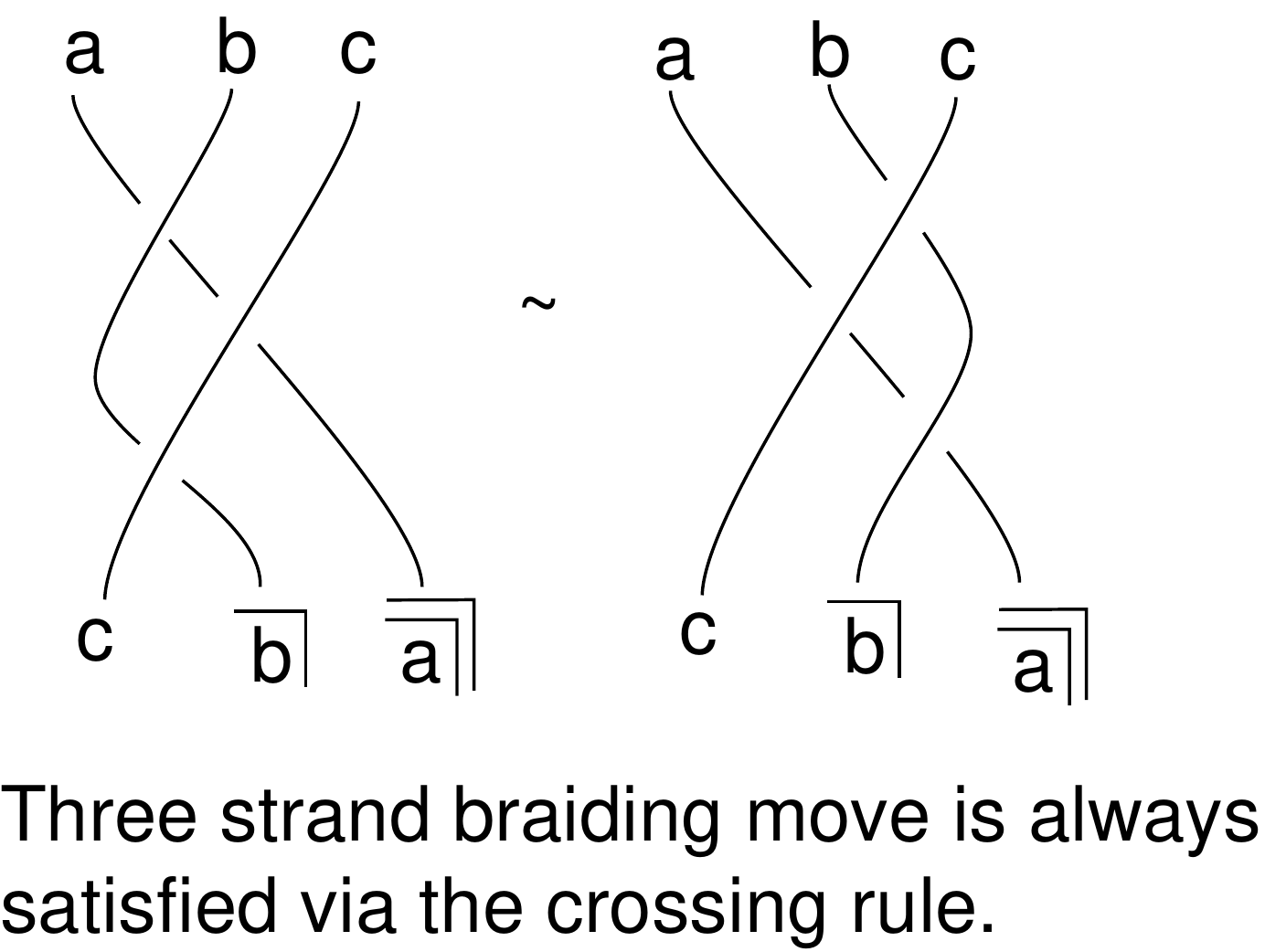}}
\caption{Fundamental Braid Relation.}
\label{fig7}
\end{figure}

Furthermore, we find that the following braiding axioms are valid, as can be seen in Figure \ref{fig7}. These are the axioms for the Artin Braid Group, and thus we have shown that this method constructs a representation of the Artin Braid group \cite{bf18, bf19}. It is understood that the braid generator $\sigma_{i}$ in the Artin Braid Group $B_{n}$ of braids on $n$ strands is represented by a left-to-right crossover of the $i$-th nad $i+1$-th strands such that all the other strands remain vertical. In our drawings in this paper we have indicated
examples for $n=2,3,4$ and we leave it to the reader to make further illustrations for other numbers of strands.\\

\begin{enumerate}
    \item $\sigma_i \sigma_j = \sigma_j \sigma_i$ for $|i - j| > 1$.\\[-6pt]
    \item $\sigma_i \sigma_j \sigma_i = \sigma_j \sigma_i \sigma_j$ for $|i - j| = 1$.
\end{enumerate}

%{It is quite useful to include the ``virtual" crossing, corresponding to the simple permutation  $\neg(a,b) = (b,a)$. The two braids in Figure \ref{fig7} \textcolor{red}{CHANGE FIGURE!!} forms above are generators for the Dihedral Group, so each Dihedral operation described in Section \ref{interpret}) can be represented in braided form.\\}

\begin{figure}[htbp]
\centerline{ \includegraphics[width=1.5 in]{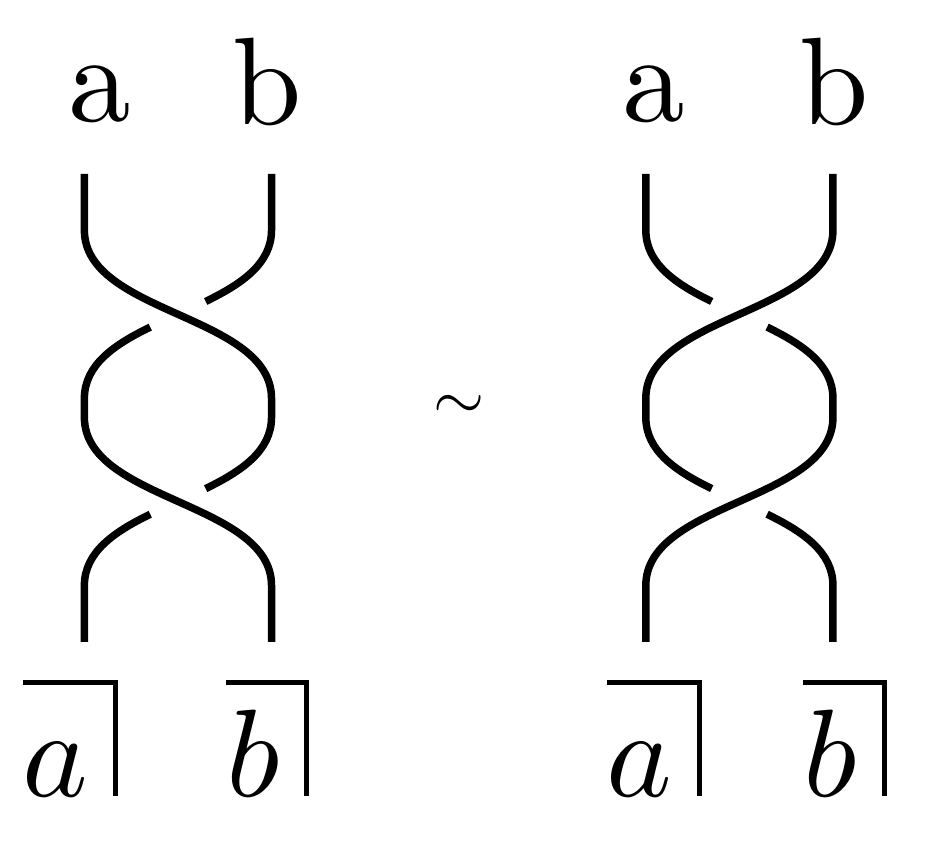}}
\caption{Doubly crossed braids of opposite polarity are equivalent.}
\label{fig8}
\end{figure}

\begin{figure}[htbp]
\centerline{ \includegraphics[width=2.8 in]{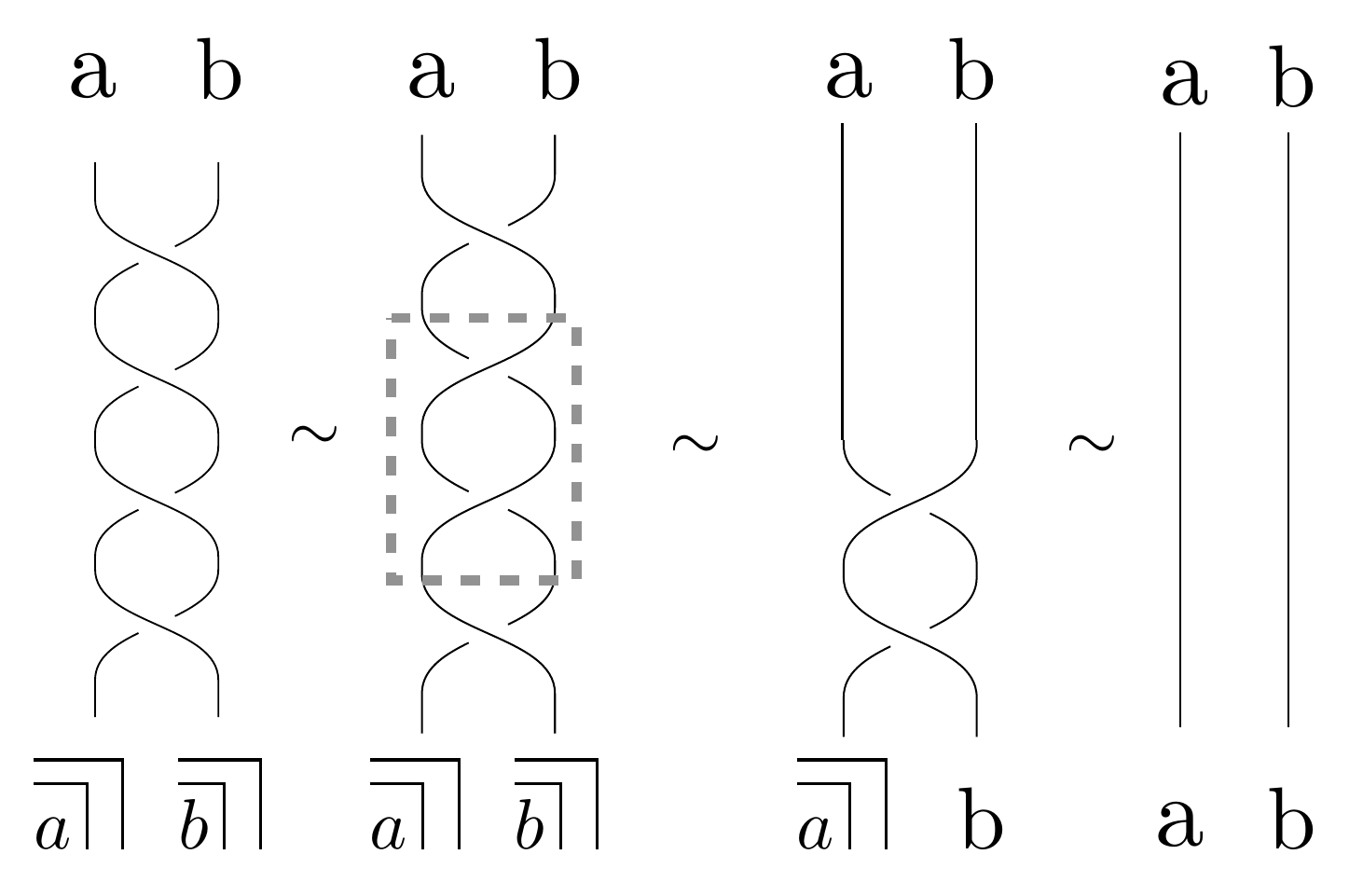}}
\caption{A quadruple braid is equivalent to the Identity Braid.}
\label{fig9}
\end{figure}

Importantly, we  see that these braided forms also satisfy additional laws that are not characteristic of the braid group in general, but are specific to the particular form of representation we have adopted. In Figure \ref{fig8}, we find that doubly crossed braids of opposite polarity are equivalent, which corresponds to the equality  $$\cross{\cross{(a,b)}_i}_i = \cross{\cross{(a,b)}_i^3}_i^3\ .$$ 

Note this is a special property of this representation of the braid group. It corresponds to the fact that $\small \cross{\quad}_i$ has order four.

In Figure \ref{fig9}, we prove that a quadruple braiding is equivalent to the identity braid. Note that in the second braid we have applied the relation in Figure \ref{fig8} within the dashed box,  replacing two twists with reverse twists.   

\begin{figure}[htbp]
\centerline{ \includegraphics[width=2.8 in]{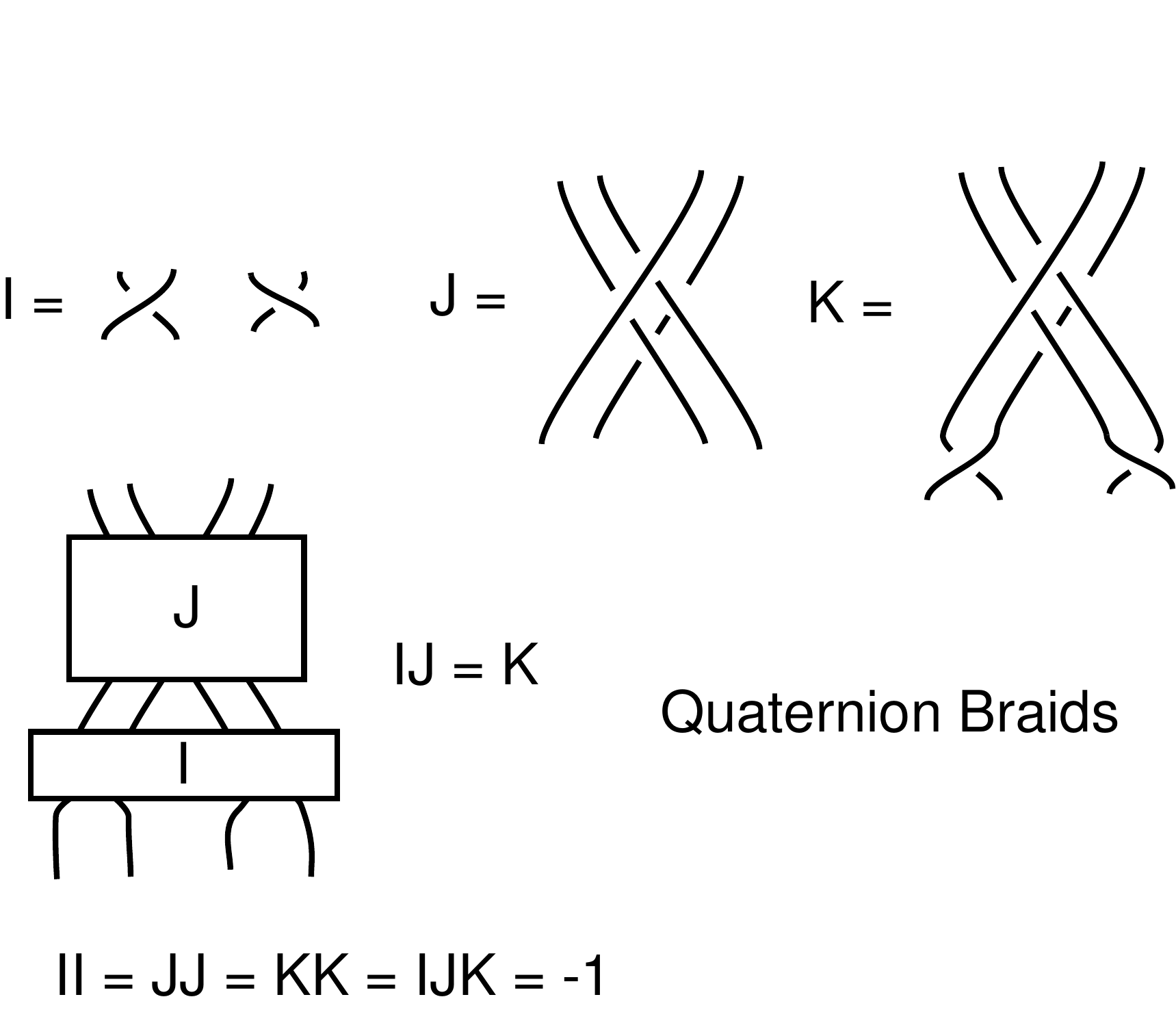}}
\caption{Braid Represention of Quaternions.}
\label{braidquaternions}
\end{figure}

The upshot of this discussion about the order four nature of the braiding generators shows that in our represenation of the braid group we can add the extra relations:
$$\sigma_{i}^2 = -1$$ where it is understood that $-1$ is an element of order two in the group that commutes with all elements of the group.
This means that we have represented the Artin braid group $B_{n}$ (whose generators have infinite order) in the group $SP_{n}$ with presentation as shown below:
\begin{enumerate}
   \item $\sigma_{i}^2 = -1,$
    \item $\sigma_i \sigma_j = \sigma_j \sigma_i$ for $|i - j| > 1$,
    \item $\sigma_i \sigma_j \sigma_i = \sigma_j \sigma_i \sigma_j$ for $|i - j| = 1$.
\end{enumerate}
The group $SP_{n}$ is a finite group, sometimes called the group of {\it signed permutations} as it can be faithfully represented by permutation matrices with entries that are
either $+1$ or $-1.$ We point out via Figure \ref{braidquaternions} that our representation can be used to give a braided interpretation for the quaternions. In this figure we illustrate
four-stranded braids $I,J,K$ so that, in the group $SP_{4}$ we have $I^2 = J^2 = K^2 = IJK = -1.$ There is more to say about this relationship with Laws of Form, BF Calculus, braiding and quaternions that we shall take up in separate publication. We have barely touched on the details of the extended algebras, but have shown enough to suggest deep relations with iterative and braided forms (see for example \cite{bf61, bf17}).  \\
%\section*{Acknowledgment}i

\end{document}